\documentclass[12pt,draft]{article}

\usepackage[cp1251]{inputenc}
\usepackage[T2A]{fontenc}
\usepackage[english,ukrainian]{babel}
\usepackage{amsmath,amsfonts,amssymb}
\usepackage{geometry}
\usepackage{cite}
\begin{document}

\begin{flushleft}

\textbf{Anna Anop, Tetiana Kasirenko, and Aleksandr Murach}\\
\small(Institute of Mathematics, National Academy of Sciences of Ukraine, Kyiv)

\medskip

\large\textbf{NONREGULAR ELLIPTIC BOUNDARY-VALUE PROBLEMS AND H\"ORMANDER SPACES}\normalsize

\medskip\medskip

\textbf{Анна Аноп, Тетяна Касіренко і Олександр Мурач}\\
\small(Інститут математики НАН України, Київ)

\medskip

\large\textbf{НЕРЕГУЛЯРНІ ЕЛІПТИЧНІ КРАЙОВІ ЗАДАЧІ ТА ПРОСТОРИ ХЕРМАНДЕРА}\normalsize

\end{flushleft}


\medskip

\noindent We investigate nonregular elliptic problems with boundary conditions of higher orders. We prove that these problems are Fredholm on appropriate pairs of inner product H\"ormander spaces that form a two-sided refined Sobolev scale. We also prove a theorem on the regularity of generalized solutions to the problems in these spaces.

\medskip

\noindent Досліджено нерегулярні еліптичні задачі з крайовими операторами вищих порядків. Доведено, що ці задачі є нетеровими у відповідних парах гільбертових просторів Хермандера, які утворюють двобічну уточнену соболєвську шкалу. Доведено також теорему про регулярність узагальнених розв'язків досліджуваних задач у цих просторах.

\bigskip

\noindent\textbf{1. Вступ.} Ця робота присвячена дослідженню в класах просторів Хермандера еліптичних задач, для яких максимум порядків крайових операторів більший за порядок еліптичного рівняння, або рівний йому. Такі задачі є нерегулярними еліптичними; для них не виконується класична формула Гріна, що ускладнює їх дослідження. Змістовні приклади цих задач зустрічаються в акустиці, гідродинаміці, теорії випадкових процесів \cite{Venttsel59, Krasil'nikov61, VeshevKouzov77}.

Зазначені задачі досить повно досліджено у двобічній шкалі просторів Соболєва, модифікованих за Ройтбергом (див. монографії \cite{Roitberg96} (розд. 4, 7) і \cite{KozlovMazyaRossmann97} (п. 4.1)). Доведено теореми про нетеровість цих задач і регулярність їх розв'язків у просторах Соболєва\,-\,Ройтберга. Останні збігаються з соболєвськими просторами, якщо їх порядок регулярності є достатньо великим числом; у противному разі вони містять елементи, які не є розподілами. Втім, соболєвська шкала, градуйована за допомогою числового показника регулярності, є занадто грубою для низки задач теорії диференціальних рівнянь і математичного аналізу. На це вказував Л.~Хермандер \cite{Hermander63, Hermander83} ще у 1963~р., який увів і дослідив широкі класи нормованих просторів, для яких показником регулярності служить досить загальна функція, та застосував їх до дослідження рівнянь з частинними похідними. В~останні двадцять років простори Хермандера та їх різні узагальнення широко застосовуються у різних розділах математики \cite{Jacob010205, MikhailetsMurach14, NicolaRodino10, Paneah00, Stepanets05, Triebel01}.

Недавно В.~А.~Михайлець і О.~О.~Мурач \cite{MikhailetsMurach05UMJ5, MikhailetsMurach06UMJ3, MikhailetsMurach06UMJ11, MikhailetsMurach06UMB4, MikhailetsMurach07UMJ5, MikhailetsMurach08UMJ4} побудували загальну теорію розв'язності еліптичних крайових задач у гільбертових просторах Хермандера $H^{s,\varphi}$, які утворюють уточнену соболєвську шкалу (їх результати підсумовано в \cite{MikhailetsMurach12BJMA2, MikhailetsMurach14}). Показниками регулярності для цих просторів служать дійсне число $s$ і функція $\varphi:(0,\infty)\to(0,\infty)$, повільно змінна на нескінченності за Караматою. Функціональний параметр $\varphi$ уточнює основну регулярність $s$. У~випадку $\varphi(\cdot)\equiv1$ простір $H^{s,\varphi}$ стає гільбертовим простором Соболєва $H^{s}$ порядку $s$. Зовсім недавно ця теорія була доповнена у статтях \cite{AnopMurach14MFAT2, AnopMurach14UMJ7, AnopKasirenko16MFAT4}, де застосовано більш широкі класи гільбертових просторів Хермандера \cite{MikhailetsMurach13UMJ3, MikhailetsMurach15ResMath1}. Відмітимо також роботу \cite{DenkFaierman17arxiv}, в якій досліджено еліптичні з параметром крайові задачі у нормованих просторах з функціональним показником регулярності.

Проте, еліптичні задачі, яким присвячено цю роботу, не були охоплені згаданою теорією. Мета нашої роботи~--- довести теореми про нетеровість досліджуваних задач і регулярність їх узагальнених розв'язків у двобічній уточненій соболєвській шкалі. У~випадку, коли числовий показник $s>m+1/2$, де $m$~--- максимум порядків крайових операторів, відповідні версії цих теорем доведено нами в \cite{KasirenkoMurach17UMJ11}. Випадок $s\leq m+1/2$ істотно більш складний для дослідження, оскільки у ньому ліві частини крайових умов не можна коректно означити на класі $H^{s,\varphi}(\Omega)$ ров'язків еліптичного рівняння, заданого в обмеженій евклідовій області $\Omega$ з гладкою межею. На відміну від згаданих монографій \cite{Roitberg96, KozlovMazyaRossmann97}, ми дотримуємося в ідейному плані підходу Ж.-Л.~Ліонса і Е.~Мадженеса \cite{LionsMagenes62, LionsMagenes63, LionsMagenes71}, розробленому для регулярних еліптичних крайових задач у двобічній соболєвській шкалі. Ми обмежуємося розглядом розв'язків $u\in H^{s,\varphi}(\Omega)$ еліптичного рівняння, права частина якого достатньо регулярна~--- належить простору Соболєва $H^{\lambda}(\Omega)$, де $\lambda>m+1/2-2q$, а $2q$~--- порядок цього рівняння. Як буде показано, крайові умови допускають коректне означення на класі усіх таких розв'язків, а відповідна еліптична крайова задача володіє властивостями, подібними до її властивостей у випадку $s>m+1/2$. Наскільки нам відомо, отримані у цій роботі результати є новими і для соболєвських просторів.

Зауважимо, що у випадку, коли порядки крайових умов менші за порядок еліптичного рівняння, різні версії теорем Ліонса\,--\,Мадженеса про нетеровість еліптичних крайових задач доведено в \cite{Roitberg68, KostarchukRoitberg73, Murach09MFAT2} для просторів Соболєва і в \cite{MikhailetsMurach11Dop4, MurachChepurukhina15UMJ} для уточненої соболєвської шкали (частина цих результатів викладена у монографії \cite{MikhailetsMurach14} (пп. 4.4, 4.5)).

\textbf{2. Постановка задачі.} Нехай $\Omega$~--- довільна обмежена область в $\mathbb{R}^{n}$, де ціле $n\geq2$. Припускаємо, що її межа $\Gamma$ є нескінченно гладким компактним многовидом вимірності $n-1$, причому $C^\infty$-структура на $\Gamma$ породжена простором~$\mathbb{R}^{n}$.

Розглянемо в області $\Omega$ таку крайову задачу:
\begin{gather}\label{1f1}
Au=f\quad\mbox{в}\quad\Omega,\\
B_{j}u=g_{j}\quad\mbox{на}\quad\Gamma,
\quad j=1,...,q. \label{1f2}
\end{gather}
Тут $A:=A(x,D)$~--- лінійний диференціальний оператор на $\overline{\Omega}:=\Omega\cup\Gamma$ довільного парного порядку $2q\geq2$, а кожне $B_{j}:=B_{j}(x,D)$~--- крайовий лінійний диференціальний оператор на $\Gamma$ довільного порядку $m_{j}\geq0$. Усі коефіцієнти цих диференціальних операторів є нескінченно гладкими комплекснозначними функціями, заданими на $\overline{\Omega}$ і $\Gamma$ відповідно. Узагалі в роботі розподіли та функції вважаємо комплекснозначними і тому розглядаємо комплексні функціональні простори.

Припускаємо, що крайова задача \eqref{1f1}, \eqref{1f2} еліптична в області $\Omega$, тобто диференціальний оператор $A$ є правильно еліптичним на $\overline{\Omega}$, а набір $B:=(B_{1},\ldots,B_q)$  крайових диференціальних операторів задовольняє умову Лопатинського щодо $A$ на $\Gamma$ (див., наприклад, огляд \cite{Agranovich97} (п.~1.2) або довідник \cite{FunctionalAnalysis72} (розд.~III, \S~6, пп. 1, 2)). Окрім того, припускаємо, що
$$
m:=\max\{m_{1},\ldots,m_{q}\}\geq2q.
$$
Отже, еліптична крайова задача \eqref{1f1}, \eqref{1f2} є нерегулярною. Задля більшої лаконічності формул покладемо $r:=m+1$.

Пов'яжемо із цією задачею лінійне відображення
\begin{equation}\label{1f3}
u\mapsto(Au,Bu)=(Au,B_{1}u,\ldots,B_{q}u),\quad\mbox{де}\quad u\in C^{\infty}(\overline{\Omega}).
\end{equation}
Будемо досліджувати властивості продовження за неперервністю цього відображення у підходящих парах функціональних просторів Хермандера.

Для опису області значень цього продовження нам знадобиться така спеціальна формула Гріна \cite{KozlovMazyaRossmann97} (формула (4.1.10)):
\begin{gather*}
(Au,v)_{\Omega}+\sum_{j=1}^{r-2q}(D_{\nu}^{j-1}Au,w_{j})_{\Gamma}+
\sum_{j=1}^{q}(B_{j}u,h_{j})_{\Gamma}=\\
=(u,A^{+}v)_{\Omega}+\sum_{k=1}^{r}\biggl(D_{\nu}^{k-1}u,K_{k}v+
\sum_{j=1}^{r-2q}R_{j,k}^{+}w_{j}+
\sum_{j=1}^{q}Q_{j,k}^{+}h_{j}\biggr)_{\Gamma}
\end{gather*}
для довільних $u,v\in C^{\infty}(\overline{\Omega})$ і $w_{1},\ldots,w_{r-2q},h_{1},\ldots,h_{q}\in C^{\infty}(\Gamma)$. Тут і далі через $(\cdot,\cdot)_{\Omega}$ і $(\cdot,\cdot)_{\Gamma}$ позначено скалярні добутки у гільбертових просторах $L_{2}(\Omega)$ і $L_{2}(\Gamma)$ функцій квадратично інтегровних відповідно на $\Omega$ і $\Gamma$ за мірою Лебега, а також продовження за неперервністю цих скалярних добутків. Окрім того, $D_{\nu}:=i\partial/\partial\nu$, де $i$~--- уявна одиниця, а $\nu$~--- поле ортів внутрішніх нормалей до межі $\Gamma$. Як звичайно, $A^{+}$~--- диференціальний оператор, формально спряжений до $A$ відносно $(\cdot,\cdot)_{\Omega}$. Окрім того, всі $R_{j,k}^{+}$ і $Q_{j,k}^{+}$ є дотичними диференціальними операторами, формально спряженими відповідно до $R_{j,k}$ і $Q_{j,k}$ відносно $(\cdot,\cdot)_{\Gamma}$. Тут дотичні лінійні диференціальні оператори $R_{j,k}:=R_{j,k}(x,D_{\tau})$ і $Q_{j,k}:=Q_{j,k}(x,D_{\tau})$ узяті із зображення крайових диференціальних операторів $D_{\nu}^{j-1}A$ і $B_{j}$ у вигляді
\begin{gather*}
D_{\nu}^{j-1}A(x,D)=\sum_{k=1}^{r}R_{j,k}(x,D_{\tau})D_{\nu}^{k-1},\quad
j=1,\ldots,r-2q,\\
B_{j}(x,D)=\sum_{k=1}^{r}Q_{j,k}(x,D_{\tau})D_{\nu}^{k-1},\quad
j=1,\ldots,q.
\end{gather*}
Зауважимо, що $\mathrm{ord}\,R_{j,k}\leq 2q+j-k$ і $\mathrm{ord}\,Q_{j,k}\leq m_{j}-k+1$, причому $R_{j,k}=0$ при $k\geq2q+j+1$ і $Q_{j,k}=0$ при $k\geq m_{j}+2$. Нарешті, кожне $K_{k}:=K_{k}(x,D)$ є деяким крайовим лінійним диференціальним оператором на $\Gamma$ порядку $\mathrm{ord}\,K_{k}\leq2q-k$ з коефіцієнтами класу $C^{\infty}(\overline{\Omega})$; при цьому $K_{k}=0$, якщо $k\geq2q+1$.

Беручи до уваги спеціальну формулу Гріна, розглянемо таку крайову
задачу в області $\Omega$ з $r-q$ додатковими невідомими функціями на межі $\Gamma$:
\begin{gather}\label{1f4}
A^{+}v=\omega\quad\mbox{в}\quad\Omega,\\
K_{k}v+\sum_{j=1}^{r-2q}R_{j,k}^{+}w_{j}+
\sum_{j=1}^{q}Q_{j,k}^{+}h_{j}=\psi_{k}\quad
\mbox{на}\quad\Gamma,\quad k=1,...,r. \label{1f5}
\end{gather}
Тут функція $v$ на $\Omega$ і $r-q$ функцій $w_{1},\ldots,w_{r-2q},h_{1},\ldots,h_{q}$ на $\Gamma$ є невідомими. Ця задача називається формально спряженою до задачі \eqref{1f1}, \eqref{1f2} відносно розглянутої спеціальної формули Гріна. Як відомо \cite{KozlovMazyaRossmann97} (теорема~4.1.1), крайова задача \eqref{1f1}, \eqref{1f2} еліптична в області $\Omega$ тоді і лише тоді, коли формально спряжена задача \eqref{1f4}, \eqref{1f5} еліптична в $\Omega$ як крайова задача з додатковими невідомими функціями на межі області.

Позначимо через $N$ лінійний простір усіх розв'язків
$u\in C^{\infty}(\overline{\Omega})$ крайової задачі \eqref{1f1}, \eqref{1f2} у випадку, коли $f=0$ в $\Omega$ і кожне $g_{j}=0$ на~$\Gamma$. Окрім того, позначимо через $N_{\star}$ лінійний простір усіх розв'язків
$$
(v,w_{1},\ldots,w_{r-2q},h_{1},\ldots,h_{q})\in C^{\infty}(\overline{\Omega})\times(C^{\infty}(\Gamma))^{r-q}
$$
формально спряженої крайової задачі \eqref{1f4}, \eqref{1f5} у випадку, коли $\omega=0$ в $\Omega$ і кожне $\psi_{k}=0$ на~$\Gamma$. Оскільки обидві ці задачі еліптичні в $\Omega$, то простори $N$ і $N_{\star}$ скінченновимірні \cite{KozlovMazyaRossmann97} (наслідок~4.1.1).

\textbf{3. Простори Хермандера.} Еліптичну крайову задачу \eqref{1f1}, \eqref{1f2} досліджуємо у придатних парах гільбертових просторів Хермандера $H^{s,\varphi}$, для яких показниками регулярності (або гладкості) служать довільні число $s\in\mathbb{R}$ і функція $\varphi\in\mathcal{M}$. Ці простори утворюють уточнену соболєвську шкалу, введену і досліджену в \cite{MikhailetsMurach05UMJ5, MikhailetsMurach06UMJ3}. Тут і надалі $\mathcal {M}$~--- множина всіх вимірних за Борелем функцій $\varphi: [1,\infty)\rightarrow(0,\infty)$, які обмежені і відокремлені від нуля на кожному компакті та повільно змінюються на нескінченності за Й.~Караматою, тобто $\varphi(\lambda t)/\varphi(t)\rightarrow 1$ при $t\rightarrow\infty$ для кожного $\lambda>0$. Повільно змінні функції добре вивчені і мають різноманітні застосування \cite{Seneta76, BinghamGoldieTeugels89}. Їх характерним прикладом служить функція
$$
\varphi(t):=(\log t)^{r_{1}}(\log\log
t)^{r_{2}}\ldots(\underbrace{\log\ldots\log}_{k\;\mbox{\small разів}}
t)^{r_{k}},\quad t\gg1,
$$
де довільно вибрано ціле число $k\geq1$ і дійсні числа  $r_{1},\ldots,r_{k}$.

Нехай $s\in\mathbb{R}$ і $\varphi\in\mathcal {M}$. Наведемо означення просторів $H^{s,\varphi}$ спочатку для $\mathbb{R}^{n}$, а потім для $\Omega$ і $\Gamma$, та зазначимо деякі їх властивості, потрібні нам. При цьому будемо слідувати монографії \cite{MikhailetsMurach14} (пп.~1.3, 2.1, 3.2).

За означенням, лінійний простір $H^{s,\varphi}(\mathbb{R}^{n})$, де ціле $n\geq1$, складається з усіх розподілів $w\in\nobreak\mathcal{S}'(\mathbb{R}^{n})$ таких, що їх перетворення Фур'є $\widehat{w}$ є функцією, яка локально інтегровна на $\mathbb {R}^{n}$ за Лебегом і задовольняє умову
$$
\int\limits_{\mathbb{R}^{n}}
\langle\xi\rangle^{2s}\varphi^{2}(\langle\xi\rangle)
|\widehat{w}(\xi)|^{2}d\xi<\infty.
$$
Тут, як звичайно, $\mathcal{S}'(\mathbb{R}^{n})$~--- лінійний топологічний простір усіх повільно зростаючих розподілів на $\mathbb{R}^{n}$, а $\langle\xi\rangle:=(1+|\xi|^{2})^{1/2}$. У просторі $H^{s,\varphi}(\mathbb {R}^{n})$ уведено скалярний добуток розподілів $w_{1}$ і $w_{2}$ за формулою
$$
(w_{1},w_{2})_{H^{s,\varphi}(\mathbb {R}^{n})}:= \int\limits_{\mathbb{R}^{n}}
\langle\xi\rangle^{2s}\varphi^{2}(\langle\xi\rangle)
\widehat{w_{1}}(\xi)\overline{\widehat{w_{2}}(\xi)}d\xi.
$$
Він породжує норму
$$
\|w\|_{H^{s,\varphi}(\mathbb{R}^{n})}:=
(w,w)_{H^{s,\varphi}(\mathbb{R}^{n})}^{1/2}.
$$

Простір $H^{s,\varphi}(\mathbb{R}^{n})$~--- ізотропний гільбертів випадок простору $\mathcal{B}_{p,k}$, введеного і дослідженого Л.~Хермандером \cite{Hermander63} (п.~2.2). А саме, $H^{s,\varphi}(\mathbb{R}^{n})=\mathcal{B}_{2,k}$, якщо $k(\xi)=\langle\xi\rangle^{s}\varphi(\langle\xi\rangle)$ для довільного $\xi\in\mathbb{R}^{n}$.

У важливому окремому випадку, коли $\varphi(\cdot)\equiv1$, простір $H^{s,\varphi}(\mathbb {R}^{n})$ стає гільбертовим простором Соболєва $H^{s}(\mathbb {R}^{n})$ порядку $s$. У загальній ситуації виконуються неперервні та щільні вкладення
\begin{equation}\label{1f7}
H^{s+\varepsilon}(\mathbb{R}^{n})\hookrightarrow H^{s,\varphi}(\mathbb{R}^{n})\hookrightarrow H^{s-\varepsilon}(\mathbb{R}^{n})
\quad\mbox{для кожного}\quad\varepsilon>0.
\end{equation}
З них випливає, що у класі функціональних просторів $\{H^{s,\varphi}(\mathbb{R}^{n}):s\in\mathbb{R},\varphi\in\mathcal{M}\}$
чи\-сло\-вий параметр $s$ задає основну регулярність (або гладкість) розподілів, а функціональний параметр $\varphi$ задає додаткову регулярність, яка уточнює основну. Тому цей клас природно називати уточненою соболєвською шкалою на $\mathbb{R}^{n}$.

Її аналоги для евклідової області $\Omega$ і замкненого
компактного многовиду $\Gamma$ вводяться у стандартний спосіб.
Наведемо відповідні означення.

За означенням, лінійний простір $H^{s,\varphi}(\Omega)$ складається зі звужень в область $\Omega$ усіх розподілів $w\in H^{s,\varphi}(\mathbb{R}^{n})$.
Норма у ньому означена за формулою
$$
\|u\|_{H^{s,\varphi}(\Omega)}:=
\inf\bigl\{\,\|w\|_{H^{s,\varphi}(\mathbb{R}^{n})}:\,
w\in
H^{s,\varphi}(\mathbb{R}^{n}),\;w=u\;\,\mbox{в}\;\,\Omega\,\bigr\},
$$
де $u\in H^{s,\varphi}(\Omega)$. Цей простір гільбертів і сепарабельний відносно вказаної норми та неперервно вкладений у топологічний простір $\mathcal{D}'(\Omega)$ усіх розподілів в $\Omega$. Множина $C^{\infty}(\overline{\Omega})$ щільна у просторі $H^{s,\varphi}(\Omega)$. Він є окремим випадком гільбертових просторів, уведених і досліджених Л.~Р.~Волевичим і Б.~П.~Панеяхом \cite[\S~2]{VolevichPaneah65}.

Коротко кажучи, простір $H^{s,\varphi}(\Gamma)$ складається з усіх розподілів на $\Gamma$, які в локальних координатах дають елементи простору $H^{s,\varphi}(\mathbb{R}^{n-1})$. Дамо докладне означення. Нехай довільним чином вибрано скінченний атлас із $C^{\infty}$-структури на многовиді $\Gamma$, утворений локальними картами $\pi_j:
\mathbb{R}^{n-1}\leftrightarrow\Gamma_{j}$, де $j=1,\ldots,p$.
Тут відкриті множини $\Gamma_{1},\ldots,\Gamma_{p}$ складають скінченне покриття многовиду $\Gamma$. Нехай, окрім того, вибрано функції $\chi_j\in\nobreak C^{\infty}(\Gamma)$, де $j=1,\ldots,p$, які утворюють розбиття одиниці на $\Gamma$, що задовольняє умову $\mathrm{supp}\,\chi_j\subset\Gamma_j$.

Тоді, за означенням, лінійний простір $H^{s,\varphi}(\Gamma)$ складається з усіх розподілів $h\in\mathcal{D}'(\Gamma)$ таких, що
$(\chi_{j}h)\circ\pi_{j}\in H^{s,\varphi}(\mathbb{R}^{n-1})$
для кожного номера $j\in\{1,\ldots,p\}$. Тут, звісно, $\mathcal{D}'(\Gamma)$~--- лінійний топологічний простір усіх розподілів на $\Gamma$, а $(\chi_{j}h)\circ\pi_{j}$ є зображенням розподілу $h$ у локальній карті $\pi_{j}$. У~просторі $H^{s,\varphi}(\Gamma)$ уведено норму за формулою
$$
\|h\|_{H^{s,\varphi}(\Gamma)}:=\biggl(\,\sum_{j=1}^{p}
\|(\chi_{j}h)\circ\pi_{j}\|_{H^{s,\varphi}(\mathbb{R}^{n-1})}^{2}
\biggr)^{1/2}.
$$
Цей простір гільбертів і сепарабельний відносно цієї норми та неперервно вкладений у $\mathcal{D}'(\Gamma)$. Важливо, що простір $H^{s,\varphi}(\Gamma)$ з точністю до еквівалентності норм не залежить від зазначеного вибору атласу і розбиття одиниці (див. \cite{MikhailetsMurach14} (теорема~2.3)). Множина $C^{\infty}(\Gamma)$ щільна у $H^{s,\varphi}(\Gamma)$.

Введені гільбертові функціональні простори утворюють уточнені соболєвські
шкали
\begin{equation}\label{1f8}
\{H^{s,\varphi}(\Omega):\,s\in\mathbb{R},\varphi\in\mathcal{M}\}
\quad\mbox{і}\quad
\{H^{s,\varphi}(\Gamma):\,s\in\mathbb{R},\varphi\in\mathcal{M}\}
\end{equation}
на $\Omega$ і $\Gamma$ відповідно. Ці шкали є двобічними за числовим параметром $s$. Вони містять двобічні гільбертові соболєвські шкали: якщо $\varphi(\cdot)\equiv1$, то $H^{s,\varphi}(\Omega)=:H^{s}(\Omega)$ і $H^{s,\varphi}(\Gamma)=:H^{s}(\Gamma)$ є простори Соболєва порядку
$s\in\mathbb{R}$. Для шкал \eqref{1f8} виконуються компактні і щільні вкладення \eqref{1f7}, якщо у формулі \eqref{1f7} замінити $\mathbb{R}^{n}$ на $\Omega$ або $\Gamma$ відповідно.

Обговоримо зв'язок між шкалами \eqref{1f8}. Нехай $s>1/2$ і $\varphi\in\mathcal{M}$; тоді відображення сліду $u\mapsto u\!\upharpoonright\!\Gamma$, де $u\in C^{\infty}(\Gamma)$, продовжується єдиним чином (за неперервністю) до обмеженого лінійного оператора $R_{\Gamma}:H^{s,\varphi}(\Omega)\rightarrow H^{s-1/2,\varphi}(\Gamma)$. Отже, для кожного розподілу $u\in H^{s,\varphi}(\Omega)$, його слід $R_{\Gamma}u$ на $\Gamma$ означений коректно. Більше того,
$$
H^{s-1/2,\varphi}(\Gamma)=\bigl\{R_{\Gamma}u:u\in H^{s,\varphi}(\Omega)\bigr\}
$$
та виконується еквівалентність норм
$$
\|h\|_{H^{s-1/2,\varphi}(\Gamma)}\asymp
\inf\,\bigl\{\|u\|_{H^{s,\varphi}(\Omega)}:\,
u\in H^{s,\varphi}(\Omega),\;\,h=R_{\Gamma}u\bigr\}
$$
на класі всіх функцій $h\in H^{s-1/2,\varphi}(\Gamma)$ (див.  \cite{MikhailetsMurach14} (теорема 3.5 і наслідок 3.1)). Проте, якщо $s<1/2$, то не можна коректно означити слід на $\Gamma$ довільного розподілу $u\in H^{s,\varphi}(\Omega)$. А саме, відображення $u\mapsto u\!\upharpoonright\!\Gamma$, де $u\in C^{\infty}(\overline{\Omega})$, не можна продовжити до неперервного лінійного оператора $R_{\Gamma}:H^{s,\varphi}(\Omega)\to\mathcal{D}'(\Gamma)$ (див.  \cite{MikhailetsMurach14} (зауваження~3.5)). Це застереження зберігає силу і для $s=1/2$ у соболєвському випадку $\varphi(\cdot)\equiv1$.

\textbf{4. Основні результати} роботи стосуються характеру розв'язності еліптичної крайової задачі \eqref{1f1}, \eqref{1f2} і регулярності її узагальнених розв'язків у просторах Хермандера, які утворюють двобічні шкали \eqref{1f8}. З огляду на зазначений вище зв'язок між цими шкалами розглянемо окремо випадки $s>m+1/2$ і $s\leq m+1/2$. У першому з них нами доведено такий результат \cite{KasirenkoMurach17UMJ11} (теорема 1).

\textbf{Твердження 1.} \it Нехай $s>m+1/2$ і $\varphi\in\mathcal{M}$. Тоді відображення \eqref{1f3} продовжується єдиним чином (за неперервністю) до обмеженого лінійного оператора
\begin{equation}\label{6f11}
(A,B):H^{s,\varphi}(\Omega)\rightarrow
H^{s-2q,\varphi}(\Omega)\oplus
\bigoplus_{j=1}^{q}H^{s-m_{j}-1/2,\varphi}(\Gamma)=:
\mathcal H_{s-2q,s,\varphi}(\Omega,\Gamma).
\end{equation}
Цей оператор нетерів. Його ядро дорівнює $N$, а область значень складається з усіх векторів
$(f,g_{1},\ldots,g_{q})\in \mathcal H_{s-2q,s,\varphi}(\Omega,\Gamma)$
таких, що
\begin{equation}\label{6f12}
\begin{gathered}
(f,v)_\Omega+\sum_{j=1}^{r-2q}(D_{\nu}^{j-1}f,w_{j})_{\Gamma}+
\sum_{j=1}^{q}(g_j,h_{j})_{\Gamma}=0 \\
\mbox{для всіх} \quad (v,w_{1},\ldots,w_{r-2q},h_{1},\ldots,h_{q})\in N_{\star}.
\end{gathered}
\end{equation}
Індекс оператора \eqref{6f11} дорівнює $\dim N-\dim N_{\star}$ та не залежить від $s$ і $\varphi$. \rm

У зв'язку з цим твердженням нагадаємо, що лінійний обмежений оператор $T:E_{1}\rightarrow E_{2}$, де $E_{1}$ і $E_{2}$~--- банахові
простори, називають нетеровим, якщо його ядро $\ker T$ і коядро $E_{2}/T(E_{1})$ скінченновимірні. Якщо цей оператор нетерів, то його область значень замкнена в $E_{2}$ і він має скінченний індекс
$$
\mathrm{ind}\,T:=\dim\ker T-\dim(E_{2}/T(E_{1})).
$$

Зі сказаного наприкінці п.~3 випливає, що умову $s>m+1/2$ у твердженні~1 не можна відкинути чи послабити. Зокрема, якщо $s\leq m+1/2$ і $\varphi(\cdot)\equiv1$, то відображення $u\mapsto B_{j}u$, де $u\in C^{\infty}(\overline{\Omega})$, не можна продовжити до неперервного лінійного оператора $B_{j}:H^{s}(\Omega)\to\mathcal{D}'(\Gamma)$ у випадку, коли $m_{j}=m$.

Для того, щоб отримати версію твердження~1 для довільного $s\leq m+1/2$, обмежимося розглядом розв'язків $u\in H^{s,\varphi}(\Omega)$ еліптичного рівняння $Au=f$, права частина якого належить до простору
\begin{equation*}
H^{m+1/2-2q+}(\Omega):=\bigcup_{\lambda>m+1/2-2q}H^{\lambda}(\Omega)=
\bigcup_{\substack{\lambda>m+1/2-2q,\\\eta\in\mathcal{M}}}
H^{\lambda,\eta}(\Omega)
\end{equation*}
(тут друга рівність виконується з огляду на вкладення \eqref{1f7}).

Нехай $s\leq m+1/2$, $\varphi\in\mathcal{M}$ і $\lambda>m+1/2-2q$. Розглянемо лінійний простір
$$
H^{s,\varphi}_{A,\lambda}(\Omega):=\bigl\{u\in H^{s,\varphi}(\Omega):
Au\in H^{\lambda}(\Omega)\bigr\},
$$
наділений нормою графіка
\begin{equation}\label{6f13}
\|u\|_{H^{s,\varphi}_{A,\lambda}(\Omega)}:=
\bigl(\,\|u\|^{2}_{H^{s,\varphi}(\Omega)}+
\|Au\|^{2}_{H^{\lambda}(\Omega)}\bigr)^{1/2}.
\end{equation}
Тут $Au$ розуміємо в сенсі теорії розподілів в області $\Omega$. У соболєвському випадку $\varphi(\cdot)\equiv1$ будемо пропускати індекс $\varphi$ у позначеннях цього та інших просторів, введених на основі просторів Хермандера $H^{s,\varphi}$.

Цей простір гільбертів відносно норми \eqref{6f13}. Справді, ця норма породжена скалярним добутком, оскільки такими є норми в правій частині рівності \eqref{6f13}. Окрім того, простір $H^{s,\varphi}_{A,\lambda}(\Omega)$ повний відносно цієї норми. Справді, якщо послідовність $(u_{k})$ фундаментальна в цьому просторі, то існують границі $u:=\lim u_{k}$ в $H^{s,\varphi}(\Omega)$ і $f:=\lim Au_{k}$ в $H^{\lambda}(\Omega)$, оскільки останні два простори повні. Диференціальний оператор $A$ неперервний у $\mathcal{D}'(\Omega)$; тому
$Au=\lim Au_{k}=f$ в $\mathcal{D}'(\Omega)$. Тут, нагадаємо, $u\in H^{s,\varphi}(\Omega)$ і $f\in H^{\lambda}(\Omega)$. Тому $u\in H^{s,\varphi}_{A,\lambda}(\Omega)$ і $\lim u_{k}=u$ у просторі $H^{s,\varphi}_{A,\lambda}(\Omega)$. Отже, цей простір повний.

\textbf{Теорема 1.} \it Нехай $s\leq m+1/2$, $\varphi\in\mathcal{M}$ і $\lambda>m+1/2-2q$. Тоді множина $C^{\infty}(\overline{\Omega})$ щільна в просторі $H^{s,\varphi}_{A,\lambda}(\Omega)$, а відображення \eqref{1f3} продовжується єдиним чином (за неперервністю) до обмеженого лінійного оператора
\begin{equation}\label{6f14}
(A,B):H^{s,\varphi}_{A,\lambda}(\Omega)\to
H^{\lambda}(\Omega)\oplus
\bigoplus_{j=1}^{q}H^{s-m_{j}-1/2,\varphi}(\Gamma)=:
\mathcal{H}_{\lambda,s,\varphi}(\Omega,\Gamma).
\end{equation}
Цей оператор нетерів. Його ядро дорівнює $N$, а область значень складається з усіх векторів $(f,g_{1},\ldots,g_{q})\in
\mathcal{H}_{\lambda,s,\varphi}(\Omega,\Gamma)$, які задовольняють умову \eqref{6f12}. Індекс оператора \eqref{6f14} дорівнює $\dim N-\dim N_{\star}$ та не залежить від $s$, $\varphi$ і~$\lambda$. \rm

Перейдемо до питання про регулярність узагальнених розв'язків крайової задачі \eqref{1f1}, \eqref{1f2} у двобічній уточненій соболєвській шкалі. Спочатку дамо означення цих розв'язків. Позначимо через $\mathcal{S}'(\Omega)$ лінійний простір звужень в область $\Omega$ усіх розподілів $w\in\mathcal{S}'(\mathbb{R}^{n})$. Покладемо
\begin{equation*}
\mathcal{S}'_{A,m+1/2-2q+}(\Omega):=
\bigl\{u\in\mathcal{S}'(\Omega):Au\in H^{m+1/2-2q+}(\Omega)\bigr\}.
\end{equation*}
Нехай розподіл $u\in\mathcal{S}'_{A,m+1/2-2q+}(\Omega)$; тоді $u\in H^{s}_{A,\lambda}(\Omega)$ для деяких чисел $s\leq m+1/2$ і $\lambda>m+1/2-2q$.  Розподіл $u$ називаємо (сильним) узагальненим розв'язком крайової задачі \eqref{1f1}, \eqref{1f2} з правою частиною
\begin{equation*}
(f,g_{1},\ldots,g_{q})\in\mathcal{S}'(\Omega)\times(\mathcal{D}'(\Gamma))^{q},
\end{equation*}
якщо $(A,B)u=(f,g_{1},\ldots,g_{q})$, де $(A,B)$~--- оператор \eqref{6f14}. Звісно, це означення не залежить від вибору чисел $s$ і~$\lambda$.

Нехай $V$~--- відкрита множина в $\mathbb{R}^{n}$, яка задовольняє умову $\Omega_0:=\Omega\cap V\neq\varnothing$. Покладемо $\Gamma_{0}:=\Gamma\cap V$ (можливий випадок, коли $\Gamma_{0}=\varnothing$). Уведемо для множин $\Omega_0$ і $\Gamma_{0}$ локальні аналоги просторів $H^{\sigma,\varphi}(\Omega)$ і $H^{\sigma,\varphi}(\Gamma)$, де $\sigma\in\mathbb{R}$ і $\varphi\in\mathcal{M}$. За означенням, лінійний простір $H^{\sigma,\varphi}_{\mathrm{loc}}(\Omega_{0},\Gamma_{0})$ cкладається з усіх розподілів $u\in\mathcal{S}'(\Omega)$ таких, що $\chi u\in\nobreak H^{\sigma,\varphi}(\Omega)$ для довільної функції $\chi\in C^{\infty}(\overline{\Omega})$, носій якої задовольняє умову $\mathrm{supp}\,\chi\subset\Omega_0\cup\Gamma_{0}$. Аналогічно, лінійний простір $H^{\sigma,\varphi}_{\mathrm{loc}}(\Gamma_{0})$ cкладається, за означенням, з усіх розподілів $h\in\mathcal{D}'(\Gamma)$ таких, що $\chi h\in H^{\sigma,\varphi}(\Gamma)$ для довільної функції $\chi\in C^{\infty}(\Gamma)$, носій якої задовольняє умову $\mathrm{supp}\,\chi\subset\Gamma_{0}$.

\textbf{Теорема 2.} \it Нехай $s\in\mathbb{R}$ і $\varphi\in\mathcal{M}$. Припустимо, що розподіл $u\in\mathcal{S}'_{A,m+1/2-2q+}(\Omega)$ є узагальненим розв'язком еліптичної крайової задачі \eqref{1f1}, \eqref{1f2}, праві частини якої задовольняють умови $f\in H^{s-2q,\varphi}_{\mathrm{loc}}(\Omega_{0},\Gamma_{0})$ і $g_{j}\in H^{s-m_{j}-1/2,\varphi}_{\mathrm{loc}}(\Gamma_{0})$ для кожного $j\in\{1,\ldots,q\}$. Тоді $u\in H^{s,\varphi}_{\mathrm{loc}}(\Omega_{0},\Gamma_{0})$. \rm

Як бачимо, уточнена регулярність $\varphi\in\mathcal{M}$ правих частин досліджуваної задачі успадковується її узагальненим розв'язком. Відмітимо важливі окремі випадки теореми~2. Якщо $\Omega_{0}=\Omega$ і $\Gamma_{0}=\Gamma$, то простори $H^{\sigma,\varphi}_{\mathrm{loc}}(\Omega_{0},\Gamma_{0})$ і
$H^{\sigma,\varphi}_{\mathrm{loc}}(\Gamma_{0})$ збігаються з просторами $H^{\sigma,\varphi}(\Omega)$ і $H^{\sigma,\varphi}(\Gamma)$ відповідно. У цьому випадку теорема~2 стверджує, що регулярність узагальненого розв'язку $u$ підвищується глобально, тобто в усій області $\Omega$ аж до її межі~$\Gamma$. Якщо $\Gamma_{0}=\varnothing$, то регулярність розв'язку $u$ підвищується в околах усіх точок $x\in\Omega_{0}$ за умови, що ці околи не перетинають межу підобласті $\Omega_{0}$. У випадку $\Gamma_{0}=\varnothing$ простір $H^{\sigma,\varphi}_{\mathrm{loc}}(\Gamma_{0})$ збігається з  $\mathcal{D}'(\Gamma)$ і тому умова на $g_{j}$ у теоремі~2 стає тривіальною. У цьому випадку висновок теореми~2 випливає з \cite{MikhailetsMurach14} (теорема 4.19). Отже,
$$
\mathcal{S}'_{A,m+1/2-2q+}(\Omega)\subset\bigcup_{\sigma>m+1/2}
H^{\sigma}_{\mathrm{loc}}(\Omega,\varnothing).
$$

Теореми 1 і 2 є новими навіть у соболєвському випадку, коли $\varphi(\cdot)\equiv1$.

\textbf{5. Доведення основних результатів.} Обґрунтуємо основні результати роботи~--- теореми 1 і 2.

\textbf{\textit{Доведення теореми} 1.} Спочатку обґрунтуємо її у соболєвському випадку, коли $\varphi(\cdot)\equiv1$ і ціле $s<2q$. У цьому випадку щільність множини $C^{\infty}(\overline{\Omega})$ у просторі $H^{s}_{A,\lambda}(\Omega)$, випливає з теорем 4.25(і) та 4.26 з монографії \cite{MikhailetsMurach14} (див. також \cite{Murach09MFAT2}, теореми 1(i) та 2). Справді, за другою з них простір $H^{\lambda}(\Omega)$ задовольняє умову $\mathrm{I}_{s-2q}$, сформульовану в п.~4.4.2 цієї монографії. Тому за першою з цих теорем множина
$$
\bigl\{u\in C^{\infty}(\overline{\Omega}):Au\in H^{\lambda}(\Omega)\bigr\}=C^{\infty}(\overline{\Omega})
$$
щільна в просторі $H^{s}_{A,\lambda}(\Omega)$.

Для доведення решти теореми~1 у розглянутому випадку скористаємося таким результатом про розв'язність досліджуваної задачі \eqref{1f1}, \eqref{1f2} у просторах Соболєва\,--\,Ройтберга: відображення \eqref{1f3} продовжується єдиним чином (за неперервністю) до нетеревого обмеженого лінійного оператора
\begin{equation}\label{6f26-m+1}
(A,B):H^{s,(r)}(\Omega)\to H^{s-2q,(r-2q)}(\Omega)\oplus
\bigoplus_{j=1}^{q}H^{s-m_{j}-1/2}(\Gamma)=:
\mathcal{H}_{s-2q,s}^{(r-2q)}(\Omega,\Gamma).
\end{equation}

Тут через $H^{\sigma,(k)}(\Omega)$, де $\sigma\in\mathbb{R}$ і $1\leq k\in\mathbb{Z}$, позначено гільбертів функціональний простір, введений Я.~А.~Ройтбергом \cite{Roitberg64} на основі соболєвських просторів (див. також його монографію \cite{Roitberg96} (п.~2.1)). За означенням, простір Соболєва\,--\,Ройтберга $H^{\sigma,(k)}(\Omega)$, де $\sigma\notin\{1/2,\ldots,k-1/2\}$, є поповненням лінійного многовиду $C^{\infty}(\overline{\Omega})$ за гільбертовою нормою
$$
\|u\|_{H^{\sigma,(k)}(\Omega)}:=
\biggl(\|u\|_{H^{\sigma,(0)}(\Omega)}^{2}+
\sum_{j=1}^{k}\;\|(D_{\nu}^{j-1}u)\!\upharpoonright\!\Gamma\|
_{H^{\sigma-j+1/2}(\Gamma)}^{2}\biggr)^{1/2}.
$$
Тут $H^{\sigma,(0)}(\Omega):=H^{\sigma}(\Omega)$, якщо $\sigma\geq0$, та $H^{\sigma,(0)}(\Omega)$~--- поповнення $C^{\infty}(\overline{\Omega})$ за гільбертовою нормою
$$
\|u\|_{H^{\sigma,(0)}(\Omega)}:=
\sup\bigl\{\,|(u,v)_{\Omega}|:\,v\in H^{-\sigma}(\Omega),\,
\|v\|_{H^{-\sigma}(\Omega)}=1\bigr\},
$$
якщо $\sigma<0$. У випадку, коли $\sigma\in\{1/2,\ldots,k-1/2\}$, простір $H^{\sigma,(k)}(\Omega)$ означається шляхом інтерполяції з параметром $1/2$ пари гільбертових просторів $H^{\sigma\mp\varepsilon,(k)}(\Omega)$, де $0<\varepsilon<1$.

Простір $H^{\sigma,(k)}(\Omega)$, де $\sigma\notin\{1/2,\ldots,k-1/2\}$, допускає такий опис \cite{Roitberg96} (лема 2.2.1): лінійне відображення
$$
T_{k}:u\mapsto\bigl(u,u\!\upharpoonright\!\Gamma,\ldots,
(D_{\nu}^{k-1}u)\!\upharpoonright\!\Gamma\bigr),\quad\mbox{де}\quad u\in C^{\infty}(\overline{\Omega}),
$$
продовжується єдиним чином (за неперервністю) до ізометричного оператора
$$
T_{k}:\,H^{\sigma,(k)}(\Omega)\rightarrow
H^{\sigma,(0)}(\Omega)\oplus
\bigoplus_{j=1}^{k}\,H^{\sigma-j+1/2}(\Gamma)=:
\Pi_{\sigma,(k)}(\Omega,\Gamma),
$$
область значень якого складається з усіх векторів
$(u_{0},u_{1},\ldots,u_{k})\in\Pi_{\sigma,(k)}(\Omega,\Gamma)$
таких, що $u_{j}=R_{\Gamma}D_{\nu}^{j-1}u_{0}$ для кожного $j\in\{1,\ldots,k\}$, що задовольняє умову $\sigma>j-1/2$.

Зауважимо \cite{Roitberg96} (п.~2.1), що при $\sigma>k-1/2$ простори $H^{\sigma,(k)}(\Omega)$ і $H^{\sigma}(\Omega)$ рівні як поповнення лінійного многовиду $C^{\infty}(\overline{\Omega})$ за еквівалентними нормами. Окрім того, виконується неперервне вкладення $H^{\sigma,(k)}(\Omega)\hookrightarrow H^{\theta,(k)}(\Omega)$, якщо $\theta<\sigma$.

Обмеженість і нетеровість оператора \eqref{6f26-m+1} доведена Ю.~В.~Костарчуком і Я.~А.~Ройтбергом в \cite{KostarchukRoitberg73} (теорема~5) для довільного дійсного $s$. Там же показано, що $N$~--- ядро цього оператора. Згідно з \cite{Roitberg96} (теорема 4.1.3) область значень оператора \eqref{6f26-m+1} складається з усіх векторів
$(f,g_{1},\ldots,g_{q})\in\mathcal{H}_{s-2q,s}^{(r-2q)}(\Omega,\Gamma)$ таких, що
\begin{equation}\label{Roitberg-range}
\begin{gathered}
(f_{0},v)_\Omega+\sum_{j=1}^{r-2q}(f_{j},w_{j})_{\Gamma}+
\sum_{j=1}^{q}(g_j,h_{j})_{\Gamma}=0\\
\mbox{для всіх}\quad
(v,w_{1},\ldots,w_{r-2q},h_{1},\ldots,h_{q})\in \mathcal{N}_{\star},
\end{gathered}
\end{equation}
а індекс дорівнює $\dim N-\dim \mathcal{N}_{\star}$. Тут $(f_{0},f_{1},\ldots,f_{r-2q}):=T_{r-2q}f$, а $\mathcal{N}_{\star}$ є деяким скінченновимірним простором, який лежить в  $C^{\infty}(\overline{\Omega})\times(C^{\infty}(\Gamma))^{r-q}$ та не залежить від~$s$. Згідно з \cite{KozlovMazyaRossmann97} (теорема 4.1.4) можна покласти $\mathcal{N}_{\star}:=N_{\star}$ в описі \eqref{Roitberg-range} області значень оператора \eqref{6f26-m+1} для цілих $s$. Покажемо, що те саме можна зробити і у формулі індексу цього оператора. Це достатньо показати для $s=0$ з огляду на незалежність індексу від $s$. Оскільки $N_{\star}\subset
C^{\infty}(\overline{\Omega})\times(C^{\infty}(\Gamma))^{r-q}$, то $N_{\star}$ можна розглядати як скінченновимірний підпростір простору
\begin{equation*}
H^{2q}(\Omega)\oplus\bigoplus_{j=1}^{r-2q}H^{2q+j-1/2}(\Gamma)\oplus
\bigoplus_{j=1}^{q}H^{m_{j}+1/2}(\Gamma).
\end{equation*}
Останній є взаємно спряженим до простору
\begin{equation}\label{Roitberg-space-represent}
\Pi_{-2q,(r-2q)}(\Omega,\Gamma)\oplus
\bigoplus_{j=1}^{q}H^{-m_{j}-1/2}(\Gamma)
\end{equation}
відносно розширення за неперервністю скалярного добутку в $L_{2}(\Omega)\oplus(L_{2}(\Gamma))^{r-q}$. За такого розгляду простір, спряжений до $N_{\star}$, збігається з факторпростором простору \eqref{Roitberg-space-represent} за підпростором усіх векторів
$$
(f_{0},f_{1},\ldots,f_{r-2q},g_{1},\ldots,g_{q})\in
\Pi_{-2q,(r-2q)}(\Omega,\Gamma)\oplus
\bigoplus_{j=1}^{q}H^{-m_{j}-1/2}(\Gamma),
$$
які задовольняють умову \eqref{Roitberg-range}, де беремо $\mathcal{N}_{\star}:=N_{\star}$. Звідси, оскільки оператор $T_{r-2q}$ здійснює ізоморфізм простору $H^{-2q,(r-2q)}(\Omega)$ на простір $\Pi_{-2q,(r-2q)}(\Omega,\Gamma)$, випливає, що ковимірність області значень оператора \eqref{6f26-m+1}, де $s=0$, дорівнює вимірності цього факторпростору, тобто становить $\dim N_{\star}$. Отже, індекс цього оператора дорівнює $\dim N-\dim N_{\star}$.

Покладемо
$$
H^{s,(r)}_{A,\lambda}(\Omega):=\bigl\{u\in H^{s,(r)}(\Omega):
Au\in H^{\lambda}(\Omega)\bigr\}.
$$
Тут для кожного $u\in H^{s,(r)}(\Omega)$ елемент $Au\in H^{s-2q,(r-2q)}(\Omega)$ означено за допомогою оператора \eqref{6f26-m+1}. Для цього елемента умова $Au\in H^{\lambda}(\Omega)$ має сенс, оскільки, як зазначалося вище,
\begin{equation}\label{R-S-embedding}
H^{\lambda}(\Omega)=H^{\lambda,(r-2q)}(\Omega)\hookrightarrow H^{s-2q,(r-2q)}(\Omega),
\end{equation}
причому вкладення неперервне. Наділимо лінійний простір $H^{s,(r)}_{A,\lambda}(\Omega)$ нормою графіка
\begin{equation}\label{6f13a}
\|u\|_{H^{s,(r)}_{A,\lambda}(\Omega)}:=
\bigl(\,\|u\|^{2}_{H^{s,(r)}(\Omega)}+
\|Au\|^{2}_{H^{\lambda}(\Omega)}\bigr)^{1/2}.
\end{equation}

Цей простір гільбертів відносно норми \eqref{6f13a}. Справді,
ця норма, звісно, породжена деяким скалярним добутком. Окрім того, простір $H^{s,(r)}_{A,\lambda}(\Omega)$ повний відносно неї. Дійсно, якщо послідовність $(u_{k})$ фундаментальна в цьому просторі, то існують границі $u:=\nobreak\lim u_{k}$ в $H^{s,(r)}(\Omega)$ і $f:=\lim Au_{k}$ в $H^{\lambda}(\Omega)$, оскільки останні два простори повні. З~першої границі випливає, що $Au=\lim Au_{k}$ в $H^{s-2q,(r-2q)}(\Omega)$. Звідси на підставі формули \eqref{R-S-embedding} і другої границі маємо рівність $Au=f$. Тому $u\in H^{s,(r)}_{A,\lambda}(\Omega)$ і $\lim u_{k}=u$ у просторі $H^{s,(r)}_{A,\lambda}(\Omega)$. Отже, цей простір повний.

Звуження відображення \eqref{6f26-m+1} на простір $H^{s,(r)}_{A,\lambda}(\Omega)$ є лінійним оператором
\begin{equation}\label{Oper-Roitberg}
(A,B):H^{s,(r)}_{A,\lambda}(\Omega)\to H^{\lambda}(\Omega)\oplus
\bigoplus_{j=1}^{q}H^{s-m_{j}-1/2}(\Gamma)=:
\mathcal{H}_{\lambda,s}(\Omega,\Gamma).
\end{equation}
Із вказаних вище властивостей оператора \eqref{6f26-m+1} безпосередньо випливає, що оператор \eqref{Oper-Roitberg} обмежений, його ядро дорівнює $N$, а область значень
\begin{equation}\label{proof1-a}
(A,B)\bigl(H^{s,(r)}_{A,\lambda}(\Omega)\bigr)=
\mathcal{H}_{\lambda,s}(\Omega,\Gamma)\cap
(A,B)\bigl(H^{s,(r)}(\Omega)\bigr).
\end{equation}
З цієї рівності та замкненості $(A,B)(H^{s,(r)}(\Omega))$ у просторі
$$
\mathcal{H}_{s-2q,s}^{(r-2q)}(\Omega,\Gamma)\hookleftarrow \mathcal{H}_{\lambda,s}(\Omega,\Gamma)
$$
випливає, що область значень оператора \eqref{Oper-Roitberg} замкнена у просторі $\mathcal{H}_{\lambda,s}(\Omega,\Gamma)$ і складається з усіх векторів
$(f,g_{1},\ldots,g_{q})\in\mathcal{H}_{\lambda,s}(\Omega,\Gamma)$, які задовольняють умову \eqref{Roitberg-range}, де $(f_{0},f_{1},\ldots,f_{r-2q}):=T_{r-2q}f$ і $\mathcal{N}_{\star}:=N_{\star}$. Оскільки $f\in H^{\lambda}(\Omega)$ і $\lambda>m+1/2-2q$, то з означення оператора $T_{r-2q}$ випливає, що $f_{0}=f$ і $f_{j}=R_{\Gamma}D_{\nu}^{j-1}f$ для кожного
$j\in\{1,\ldots,r-2q\}$. Тому умова \eqref{Roitberg-range} набирає вигляду \eqref{6f12}.

Для обґрунтування нетеровості оператора \eqref{Oper-Roitberg} залишається показати, що його область значень має скінченну ковимірність. Нагадаємо, що ковимірність області значень оператора \eqref{6f26-m+1} дорівнює $\dim N_{\star}<\infty$. Окрім того, множина
$C^{\infty}(\overline{\Omega})\times(C^{\infty}(\Gamma))^{q}$ щільна в просторі $\mathcal{H}_{s-2q,s}^{(r-2q)}(\Omega,\Gamma)$. Тому, за лемою Гохберга\,--\,Крейна \cite{HohbergKrein57} (лема 2.1), існує скінченновимірний простір $N_{1}\subset C^{\infty}(\overline{\Omega})\times(C^{\infty}(\Gamma))^{q}$ такий, що
\begin{equation*}
\mathcal{H}_{s-2q,s}^{(r-2q)}(\Omega,\Gamma)=
(A,B)\bigl(H^{s,(r)}(\Omega)\bigr)\dotplus N_{1}
\end{equation*}
(як звичайно, знак $\dotplus$ служить для позначення прямої суми підпросторів). При цьому $\dim N_{1}=\dim N_{\star}$. Звуження цієї суми на простір $\mathcal{H}_{\lambda,s}(\Omega,\Gamma)$ дає на підставі формул \eqref{proof1-a} і $N_{1}\subset\mathcal{H}_{\lambda,s}(\Omega,\Gamma)$ рівність
\begin{equation}\label{proof1-b}
\mathcal{H}_{\lambda,s}(\Omega,\Gamma)=
(A,B)\bigl(H^{s,(r)}_{A,\lambda}(\Omega)\bigr)\dotplus N_{1}.
\end{equation}
Отже, ковимірність області значень оператора \eqref{Oper-Roitberg} дорівнює $\dim N_{1}=\dim N_{\star}<\infty$. Таким чином, цей оператор нетерів.

Для того, щоб завершити доведення теореми~1 у розглянутому випадку,  достатньо показати, що множина $C^{\infty}(\overline{\Omega})$ щільна у просторі $H^{s,(r)}_{A,\lambda}(\Omega)$ та норми у просторах $H^{s}_{A,\lambda}(\Omega)$ і $H^{s,(r)}_{A,\lambda}(\Omega)$ еквівалентні на цій щільній множині. Справді, тоді нетерів оператор \eqref{Oper-Roitberg} стає оператором \eqref{6f14} з формулювання цієї теореми.

Доведемо спочатку, що множина $C^{\infty}(\overline{\Omega})$ щільна в $H^{s,(r)}_{A,\lambda}(\Omega)$. Позначимо через $Q^{s,(r)}_{A,\lambda}(\Omega)$ ортогональне доповнення підпростору $N$ у гільбертовому просторі $H^{s,(r)}_{A,\lambda}(\Omega)$. Звуження нетерового оператора \eqref{Oper-Roitberg} на підпростір $Q^{s,(r)}_{A,\lambda}(\Omega)$ є ізоморфізмом
\begin{equation}\label{proof1-c}
(A,B):Q^{s,(r)}_{A,\lambda}(\Omega)\leftrightarrow
(A,B)\bigl(H^{s,(r)}_{A,\lambda}(\Omega)\bigr);
\end{equation}
тут, звісно, $(A,B)(H^{s,(r)}_{A,\lambda}(\Omega))$ трактується як підпростір простору $\mathcal{H}_{\lambda,s}(\Omega,\Gamma)$. Позначимо через $P$ оператор косого проектування простору $\mathcal{H}_{\lambda,s}(\Omega,\Gamma)$ на підпростір $(A,B)(H^{s,(r)}_{A,\lambda}(\Omega))$ паралельно підпростору $N_{1}$ з прямої суми \eqref{proof1-b}.

Подамо довільний елемент $u\in H^{s,(r)}_{A,\lambda}(\Omega)$ у вигляді $u=v+w$, де $v\in Q^{s,(r)}_{A,\lambda}(\Omega)$ і $w\in N$. Оскільки множина $C^{\infty}(\overline{\Omega})\times(C^{\infty}(\Gamma))^{q}$ щільна у просторі $\mathcal{H}_{\lambda,s}(\Omega,\Gamma)$, то для вектора
$F:=(A,B)v\in\mathcal{H}_{\lambda,s}(\Omega,\Gamma)$ існує послідовність $(F_{k})\subset C^{\infty}(\overline{\Omega})\times(C^{\infty}(\Gamma))^{q}$ така, що $F_{k}\to F$ у просторі $\mathcal{H}_{\lambda,s}(\Omega,\Gamma)$. Тоді $PF_{k}\to PF=F$ у цьому просторі, причому $(PF_{k})\subset C^{\infty}(\overline{\Omega})\times(C^{\infty}(\Gamma))^{q}$. Отже,
$$
v_{k}:=(A,B)^{-1}PF_{k}\to(A,B)^{-1}F=v\quad\mbox{в}\;\;
Q^{s,(r)}_{A,\lambda}(\Omega);
$$
тут через $(A,B)^{-1}$ позначено оператор, обернений до ізоморфізму \eqref{proof1-c}. Оскільки
$$
(A,B)v_{k}=PF_{k}\in\mathcal{H}_{\sigma-2q,\sigma}^{(r-2q)}(\Omega,\Gamma)
\quad\mbox{для кожного}\quad \sigma\in\mathbb{R},
$$
то згідно з \cite{Roitberg96} (теорема 7.1.1) виконується включення
$$
v_{k}\in\bigcap_{\sigma\in\mathbb{R}}H^{\sigma,(r)}(\Omega)=
C^{\infty}(\overline{\Omega}).
$$
Отже, $C^{\infty}(\overline{\Omega})\ni v_{k}+w\to v+w=u$ у просторі $H^{s,(r)}_{A,\lambda}(\Omega)$. З огляду на довільність елемента $u\in H^{s,(r)}_{A,\lambda}(\Omega)$ доведено щільність множини $C^{\infty}(\overline{\Omega})$ в $H^{s,(r)}_{A,\lambda}(\Omega)$.

Доведемо тепер, що норми у просторах $H^{s}_{A,\lambda}(\Omega)$ і $H^{s,(r)}_{A,\lambda}(\Omega)$ еквівалентні на $C^{\infty}(\overline{\Omega})$. Зауважимо спочатку, що
\begin{equation}\label{proof1-d-a}
\|u\|_{H^{s}_{A,\lambda}(\Omega)}\leq
\|u\|_{H^{s,(r)}_{A,\lambda}(\Omega)}\quad\mbox{для довільного}\quad
u\in C^{\infty}(\overline{\Omega}).
\end{equation}
Це випливає з означення норми у просторі $H^{s,(r)}(\Omega)$: у випадку $s\geq0$~--- безпосередньо, а у випадку $s<0$~--- з огляду на те, що
\begin{equation*}
\|u\|_{H^{s}(\Omega)}\leq\|\mathcal{O}u\|_{H^{s}(\mathbb{R}^{n})}=
\|u\|_{H^{s,(0)}(\Omega)}\quad\mbox{для довільного}\quad
u\in C^{\infty}(\overline{\Omega}).
\end{equation*}
Тут $\mathcal{O}u$ позначає продовження нулем на $\mathbb{R}^{n}$ функції $u\in C^{\infty}(\overline{\Omega})$, а рівність правильна, як зазначено в \cite[с.~52]{Roitberg96}.

Доведемо обернену оцінку до \eqref{proof1-d-a}. Для цього розглянемо еліптичну крайову задачу, яка складається з рівняння \eqref{1f1} і крайових умов
\begin{equation}\label{proof1-d}
D_{\nu}^{j-1}u=g_{j}\quad\mbox{на}\quad\Gamma,
\quad j=1,...,q.
\end{equation}
Згідно з \cite{Roitberg96} (теорема 4.1.1) відображення
\begin{equation}\label{proof1-e}
u\mapsto(Au,\mathrm{D}u):=\bigl(Au,u\!\upharpoonright\!\Gamma,\ldots,
(D_{\nu}^{q-1}u)\!\upharpoonright\!\Gamma\bigr),\quad\mbox{де}\quad
u\in C^{\infty}(\overline{\Omega}),
\end{equation}
продовжується єдиним чином (за неперервністю) до нетерового обмеженого лінійного оператора
\begin{equation}\label{proof1-f}
(A,\mathrm{D}):H^{s,(r)}(\Omega)\to H^{s-2q,(r-2q)}(\Omega)\oplus
\bigoplus_{j=1}^{q}\,H^{s-j+1/2}(\Gamma).
\end{equation}
До того ж, ядро $\mathcal{N}$ цього оператора лежить в $C^{\infty}(\overline{\Omega})$. Звуження оператора \eqref{proof1-f} на простір $H^{s,(r)}_{A,\lambda}(\Omega)$ є нетеровим обмеженим оператором
\begin{equation}\label{proof1-g}
(A,\mathrm{D}):H^{s,(r)}_{A,\lambda}(\Omega)\to H^{\lambda}(\Omega)\oplus\bigoplus_{j=1}^{q}\,H^{s-j+1/2}(\Gamma)=: \mathcal{H}^{\mathrm{D}}_{\lambda,s}(\Omega,\Gamma).
\end{equation}
Це доводиться так само як і нетеровість оператора \eqref{Oper-Roitberg}.

Крім того, згідно з \cite{MikhailetsMurach14} (теорема~4.27)   відображення \eqref{proof1-e}
продовжується єдиним чином (за неперервністю) до нетерового обмеженого лінійного оператора
\begin{equation}\label{proof1-h}
(A,\mathrm{D}):H^{s}_{A,\lambda}(\Omega)\to
\mathcal{H}^{\mathrm{D}}_{\lambda,s}(\Omega,\Gamma),
\end{equation}
причому ядро цього оператора лежить в $C^{\infty}(\overline{\Omega})$ (див. також \cite{Murach09MFAT2}, наслідок~3). Зауважимо, що зазначена теорема доведена в \cite{MikhailetsMurach14} для регулярних еліптичних крайових задач, до яких і належить задача \eqref{1f1}, \eqref{proof1-d}. Оператори \eqref{proof1-g} і \eqref{proof1-h} мають спільне ядро $\mathcal{N}\subset C^{\infty}(\overline{\Omega})$ і спільну область значень, бо вона є замиканням множини $\{(A,\mathrm{D})u:u\in C^{\infty}(\overline{\Omega})\}$ у просторі $\mathcal{H}^{\mathrm{D}}_{\lambda,s}(\Omega,\Gamma)$. Позначимо цю спільну область значень через $\mathcal{R}_{\lambda,s}(\Omega,\Gamma)$.

Нетерові оператори \eqref{proof1-g} і \eqref{proof1-h} породжують у канонічний спосіб ізоморфізми
\begin{gather*}
(A,\mathrm{D}):H^{s,(r)}_{A,\lambda}(\Omega)/\mathcal{N}\leftrightarrow
\mathcal{R}_{\lambda,s}(\Omega,\Gamma),\\
(A,\mathrm{D}):H^{s}_{A,\lambda}(\Omega)/\mathcal{N}\leftrightarrow
\mathcal{R}_{\lambda,s}(\Omega,\Gamma).
\end{gather*}
Тут, звісно, трактуємо $\mathcal{R}_{\lambda,s}(\Omega,\Gamma)$ як підпростір простору $\mathcal{H}^{\mathrm{D}}_{\lambda,s}(\Omega,\Gamma)$. Для довільної функції $u\in C^{\infty}(\overline{\Omega})$ розглянемо її клас суміжності $\widetilde{u}:=\{u+w:w\in\mathcal{N}\}$, який належить обом факторпросторам $H^{s,(r)}_{A,\lambda}(\Omega)/\mathcal{N}$ і $H^{s}_{A,\lambda}(\Omega)/\mathcal{N}$. На підставі цих ізоморфізмів маємо еквівалентність норм
\begin{equation}\label{proof1-i}
\|\widetilde{u}\|_{H^{s,(r)}_{A,\lambda}(\Omega)/\mathcal{N}}\asymp
\|(A,\mathrm{D})\widetilde{u}\|_{\mathcal{R}_{\lambda,s}(\Omega,\Gamma)}
\asymp\|\widetilde{u}\|_{H^{s}_{A,\lambda}(\Omega)/\mathcal{N}}
\end{equation}
на функціях $u\in C^{\infty}(\overline{\Omega})$. Спираючись на неї, доведену оцінку, обернену до \eqref{proof1-d-a}.

Для довільного $u\in C^{\infty}(\overline{\Omega})$ існує таке
$w\in\mathcal{N}$, що
\begin{equation}\label{proof1-k}
\|u+w\|_{H^{s,(r)}_{A,\lambda}(\Omega)}\leq
2\,\|\widetilde{u}\|_{H^{s,(r)}_{A,\lambda}(\Omega)/\mathcal{N}}.
\end{equation}
Скориставшись еквівалентністю норм на скінченновимірному просторі $\mathcal{N}$ та послідовно формулами \eqref{proof1-d-a}, \eqref{proof1-k} і \eqref{proof1-i}, отримаємо такі нерівності:
\begin{gather*}
\|u\|_{H^{s,(r)}_{A,\lambda}(\Omega)}\leq
\|u+w\|_{H^{s,(r)}_{A,\lambda}(\Omega)}+
\|w\|_{H^{s,(r)}_{A,\lambda}(\Omega)}\leq
\|u+w\|_{H^{s,(r)}_{A,\lambda}(\Omega)}+
c_{1}\|w\|_{H^{s}_{A,\lambda}(\Omega)}\leq\\
\leq\|u+w\|_{H^{s,(r)}_{A,\lambda}(\Omega)}+
c_{1}\|u+w\|_{H^{s}_{A,\lambda}(\Omega)}+
c_{1}\|u\|_{H^{s}_{A,\lambda}(\Omega)}\leq\\
\leq(1+c_{1})\|u+w\|_{H^{s,(r)}_{A,\lambda}(\Omega)}+
c_{1}\|u\|_{H^{s}_{A,\lambda}(\Omega)}\leq
2(1+c_{1})\|\widetilde{u}\|_{H^{s,(r)}_{A,\lambda}(\Omega)/\mathcal{N}}+
c_{1}\|u\|_{H^{s}_{A,\lambda}(\Omega)}\leq\\
\leq c_{2}\|\widetilde{u}\|_{H^{s}_{A,\lambda}(\Omega)/\mathcal{N}}+
c_{1}\|u\|_{H^{s}_{A,\lambda}(\Omega)}\leq
(c_{2}+c_{1})\|u\|_{H^{s}_{A,\lambda}(\Omega)};
\end{gather*}
тут $c_{1}$ і $c_{2}$~--- деякі додатні числа, які не залежать від функцій $u$ і $w$. Отже, доведено оцінку, обернену до \eqref{proof1-d-a}.

Таким чином, теорему~1 доведено у розглянутому випадку, коли  $\varphi(\cdot)\equiv1$, $s\in\mathbb{Z}$ і $s<2q$.

У загальній ситуації виведемо цю теорему з розглянутого випадку і твердження~1 за допомогою інтерполяції з функціональним параметром пар гільбертових просторів (означення цієї інтерполяції та потрібні нам її властивості  наведено, наприклад, в \cite{MikhailetsMurach14} (п.~1.1) або \cite{MikhailetsMurach08MFAT1} (п.~2)).

Виберемо ціле число $l\geq1$, яке задовольняє умови $s>-2q(l-1)$ і $\lambda<2ql$. Скористаємося нетеровим оператором \eqref{6f14} у випадку (вже розглянутому), коли $\varphi(\cdot)\equiv1$ і число $-2q(l-1)$ узято замість $s$, та нетеровим оператором \eqref{6f11} з твердження~1 у випадку, коли $\varphi(\cdot)\equiv1$ і число $\lambda+2q$ узято замість $s$. Отже, отримаємо нетерові обмежені оператори
\begin{gather}\label{6f31}
(A,B):H^{-2q(l-1)}_{A,\lambda}(\Omega)\to H^{\lambda}(\Omega)\oplus
\bigoplus_{j=1}^{q}H^{-2q(l-1)-m_{j}-1/2}(\Gamma)=
\mathcal{H}_{\lambda,-2q(l-1)}(\Omega,\Gamma),\\
(A,B):H^{\lambda+2q}(\Omega)\to H^{\lambda}(\Omega)\oplus
\bigoplus_{j=1}^{q}H^{\lambda+2q-m_{j}-1/2}(\Gamma)=
\mathcal{H}_{\lambda,\lambda+2q}(\Omega,\Gamma). \label{6f32}
\end{gather}
Вони мають спільне ядро $N$ і однаковий індекс, рівний $\dim N-\dim N_{\star}$. Окрім того, перший оператор є розширенням другого.

Покладемо $\varepsilon:=s+2q(l-1)>0$ і $\delta:=\lambda+2q-s>0$ та означимо функцію $\psi:(0,\infty)\to(0,\infty)$ за формулами $\psi(t):=t^{\varepsilon/(\varepsilon+\delta)}
\varphi(t^{1/(\varepsilon+\delta)})$, якщо $t\geq 1$, та $\psi(t):=\varphi(1)$, якщо $0<t<1$. Згідно з \cite{MikhailetsMurach14} (теорема~1.14) ця функція є інтерполяційним параметром. Застосувавши інтерполяцію з функціональним параметром $\psi$ до обмежених лінійних операторів \eqref{6f31} і \eqref{6f32}, отримаємо обмежений лінійний оператор
\begin{equation}\label{6f34}
(A,B):\bigl[H^{-2q(l-1)}_{A,\lambda}(\Omega),
H^{\lambda+2q}(\Omega)\bigr]_{\psi}
\to\bigl[\mathcal{H}_{\lambda,-2q(l-1)}(\Omega,\Gamma),
\mathcal{H}_{\lambda,\lambda+2q}(\Omega,\Gamma)\bigr]_{\psi}.
\end{equation}
Він є звуженням відображення \eqref{6f31} на інтерполяційний простір
\begin{equation}\label{6f34-dom}
[H^{-2q(l-1)}_{A,\lambda}(\Omega),H^{\lambda+2q}(\Omega)]_{\psi}.
\end{equation}
Тут і далі у доведенні через $[X_{0},X_{1}]_{\psi}$ позначено гільбертів простір, який є результатом інтерполяції з функціональним параметром $\psi$ припустимої впорядкованої пари сепарабельних гільбертових просторів $X_{0}$ і $X_{1}$ (див. означення цієї інтерполяції в \cite{MikhailetsMurach14} (п.~1.1.1)). Зауважимо, що цю пару називають припустимою, якщо виконується неперервне і щільне вкладення $X_{1}\hookrightarrow X_{0}$. Тоді також виконуються неперервні і щільні вкладення $X_{1}\hookrightarrow[X_{0},X_{1}]_{\psi}\hookrightarrow X_{0}$ (див. \cite{MikhailetsMurach14} (теорема~1.1)).
Наведені у формулі \eqref{6f34} пари гільбертових просторів є припустимими. Припустимість першої з них випливає, зокрема, із доведеної вище щільності множини $C^{\infty}(\overline{\Omega})$ у просторі $H^{-2q(l-1)}_{A,\lambda}(\Omega)$. Крім того, цей простір сепарабельний, що випливає з нетеровості оператора \eqref{6f31} і сепарабельності простору $\mathcal{H}_{\lambda,-2q(l-1)}(\Omega,\Gamma)$. Припустимість другої пари очевидна.

Множина $C^{\infty}(\overline{\Omega})$ щільна у просторі \eqref{6f34-dom}, оскільки у нього неперервно і щільно вкладено простір $H^{\lambda+2q}(\Omega)$. Тому оператор \eqref{6f34} є продовженням за неперервністю відображення \eqref{1f3}. На підставі зазначених вище властивостей нетерових операторів \eqref{6f31} і \eqref{6f32} робимо висновок згідно з теоремою про інтерполяцію нетерових операторів (див. \cite{MikhailetsMurach14} (теорема~1.5)), що  оператор \eqref{6f34} нетерів з ядром $N$, індексом $\dim N-\dim N_{\star}$ та областю значень, рівною
\begin{equation*}
\bigl[\mathcal{H}_{\lambda,-2q(l-1)}(\Omega,\Gamma),
\mathcal{H}_{\lambda,\lambda+2q}(\Omega,\Gamma)\bigr]_{\psi}\cap
(A,B)\bigl(H^{-2q(l-1)}_{A,\lambda}(\Omega)\bigr).
\end{equation*}
З останньої властивості і вже обґрунтованого опису області значень оператора \eqref{6f31} (з використанням умови \eqref{6f12}) випливає, що область значень оператора \eqref{6f34} складається з усіх векторів
\begin{equation*}
(f,g_{1},\ldots,g_{q})\in
\bigl[\mathcal{H}_{\lambda,-2q(l-1)}(\Omega,\Gamma),
\mathcal{H}_{\lambda,\lambda+2q}(\Omega,\Gamma)\bigr]_{\psi},
\end{equation*}
які задовольняють умову \eqref{6f12}.

Для того, щоб завершити доведення теореми~1, залишається показати, що простори, у яких діє нетерів оператор \eqref{6f34}, задовольняють рівності
\begin{gather}\label{6f38}
\bigl[H^{-2q(l-1)}_{A,\lambda}(\Omega),
H^{\lambda+2q}(\Omega)\bigr]_{\psi}=H^{s,\varphi}_{A,\lambda}(\Omega),\\
\bigl[\mathcal{H}_{\lambda,-2q(l-1)}(\Omega,\Gamma),
\mathcal{H}_{\lambda,\lambda+2q}(\Omega,\Gamma)\bigr]_{\psi}=
\mathcal{H}_{\lambda,s,\varphi}(\Omega,\Gamma) \label{6f38-b}
\end{gather}
з точністю до еквівалентності норм.

Рівність \eqref{6f38-b} випливає з теореми 1.5 (про інтерполяцію прямих сум просторів) і теореми 2.2 (про інтерполяцію  з функціональним параметром соболєвських просторів на $\Gamma$), наведених у \cite{MikhailetsMurach14}. А саме,
\begin{gather*}
\bigl[\mathcal{H}_{\lambda,-2q(l-1)}(\Omega,\Gamma),
\mathcal{H}_{\lambda,\lambda+2q}(\Omega,\Gamma)\bigr]_{\psi}=\\
=\bigl[H^{\lambda}(\Omega),H^{\lambda}(\Omega)\bigr]_{\psi}\oplus
\bigoplus_{j=1}^{q}\bigl[H^{-2q(l-1)-m_{j}-1/2}(\Gamma),
H^{\lambda+2q-m_{j}-1/2}(\Gamma)\bigr]_{\psi}=\\
=H^{\lambda}(\Omega)\oplus\bigoplus_{j=1}^{q}
\bigl[H^{s-m_{j}-1/2-\varepsilon}(\Gamma),
H^{s-m_{j}-1/2+\delta}(\Gamma)\bigr]_{\psi}=\\
=H^{\lambda}(\Omega)\oplus\bigoplus_{j=1}^{q}
H^{s-m_{j}-1/2,\varphi}(\Gamma)=
\mathcal{H}_{\lambda,s,\varphi}(\Omega,\Gamma)
\end{gather*}
з точністю до еквівалентності норм.

Для доведення рівності \eqref{6f38} скористаємося одним результатом про інтерполяцію підпросторів, пов'язаних з довільним обмеженим лінійним оператором, який діє в парі гільбертових просторів. Нехай $H$, $\Phi$ і $\Psi$~--- гільбертові простори, причому виконується неперервне вкладення $\Phi\hookrightarrow\Psi$. Нехай також задано обмежений лінійний оператор $T:\nobreak H\rightarrow\Psi$. Покладемо $
(H)_{T,\Phi}:=\{u\in H:\,Tu\in\Phi\}$.
Простір $(H)_{T,\Phi}$ є гільбертовим відносно норми графіка
$$
\|u\|_{(H)_{T,\Phi}}:=\bigl(\|u\|_{H}^{2}+\|Tu\|_{\Phi}^{2}\bigr)^{1/2}.
$$

\textbf{Твердження 2.} \it Нехай задано шість сепарабельних гільбертових просторів $X_{0}$, $Y_{0}$, $Z_{0}$, $X_{1}$, $Y_{1}$ і $Z_{1}$ та три лінійних відображення $T$, $R$ і $S$, що задовольняють такі сім умов:

\rm{(i)} \it пари $X=[X_{0},X_{1}]$ і $Y=[Y_{0},Y_{1}]$ припустимі;

\rm{(ii)} \it простори $Z_{0}$ і $Z_{1}$ є підпросторами деякого лінійного простору $E$;

\rm{(iii)} \it виконуються неперервні вкладення $Y_{j}\hookrightarrow Z_{j}$ при $j\in\{0,1\}$;

\rm{(iv)} \it відображення $T$ означено на $X_{0}$ і задає  обмежені оператори $T:X_{j}\rightarrow Z_{j}$ при $j\in\{0,1\}$;

\rm{(v)} \it відображення $R$ означено на $E$ і задає обмежені оператори $R:Z_{j}\rightarrow X_{j}$ при $j\in\{0,1\}$;

\rm{(vi)} \it відображення $S$ означено на $E$ і задає обмежені оператори $S:Z_{j}\rightarrow Y_{j}$ при $j\in\{0,1\}$;

\rm{(vii)} \it для кожного $\omega\in E$ виконується рівність $TR\,\omega=\omega+S\omega$.

\noindent Тоді пара просторів $[(X_{0})_{T,Y_{0}},(X_{1})_{T,Y_{1}}]$ припустима і для довільного інтерполяційного параметра $\psi\in\mathcal{B}$ виконується така рівність просторів з точністю до еквівалентності норм:
\begin{equation*}
\bigl[(X_{0})_{T,Y_{0}},(X_{1})_{T,Y_{1}}\bigr]_{\psi} =\bigl([X_{0},X_{1}]_{\psi}\bigr)_{T,[Y_{0},Y_{1}]_{\psi}}.
\end{equation*}
\rm

Аналог цього твердження був уперше встановлений Ж.-Л.~Ліонсом і Е.~Мадженесом \cite{LionsMagenes71} (теорема 14.3) для комплексної інтерполяції з числовим параметром. Для інтерполяції з функціональним параметром твердження~2 доведено в
\cite{MikhailetsMurach06UMJ11} (п.~4) (див. також \cite{MikhailetsMurach14} (п.~3.3.2)).

У твердженні~2 покладемо $X_{0}:=H^{-2q(l-1)}(\Omega)$, $X_{1}:=H^{\lambda+2q}(\Omega)$, $Y_{0}:=Y_{1}:=Z_{1}:=H^{\lambda}(\Omega)$, $Z_{0}:=E:=H^{-2ql}(\Omega)$ і $T:=A$. Тоді
\begin{equation}\label{6f39}
H_{A,\lambda}^{-2q(l-1)}(\Omega)=(X_{0})_{T,Y_{0}}\quad\mbox{і}\quad
H^{\lambda+2q}(\Omega)=(X_{1})_{T,Y_{1}}.
\end{equation}
Зауважимо, що остання рівність виконується з точністю до еквівалентності норм, оскільки $A$ є обмеженим оператором у парі просторів $H^{\lambda+2q}(\Omega)$ і $H^{\lambda}(\Omega)$ (для довільного дійсного $\lambda$). Звісно, умови (i)\,--\,(iv) твердження~2 виконуються. Побудуємо оператори $R$ і $S$, які задовольняють решту умов (v)\,--\,(vii).

Для цього скористаємося тим, що відображення $u\mapsto A^{l}A^{l+}u+u$ задає ізоморфізм
\begin{equation}\label{6f40}
A^{l}A^{l+}+I:\,H^{\sigma}_{\mathrm{D}}(\Omega)\leftrightarrow
H^{\sigma-4ql}(\Omega)\quad\mbox{для довільного}\quad\sigma\geq2ql
\end{equation}
(див., наприклад, лему~3.1 \cite{MikhailetsMurach14}, доведення якої проходить і для дійсних $\sigma\geq2ql$). Тут, як звичайно, $A^{l}$ є $l$-тою ітерацією оператора $A$, а $A^{l+}$ є формально спряженим оператором до диференціального оператора $A^{l}$ відносно скалярного добутку в $L_{2}(\Omega)$, та $I$~--- тотожний оператор. Окрім того,
$$
H^{\sigma}_{\mathrm{D}}(\Omega):=\bigl\{u\in H^{\sigma}(\Omega):
R_{\Gamma}D_{\nu}^{j-1}u=0\;\,\mbox{для кожного}\;\,j\in\{1,\ldots,2ql\}\bigr\}
$$
є підпростором простору $H^{\sigma}(\Omega)$. Оператор, обернений до \eqref{6f40}, є обмеженим лінійним оператором
\begin{equation}\label{6f41}
(A^{l}A^{l+}+I)^{-1}:\,H^{\theta}(\Omega)\rightarrow
H^{\theta+4ql}(\Omega)\quad\mbox{для довільного}\quad \theta\geq-2ql.
\end{equation}
Покладемо
$$
R:=A^{l-1}A^{l+}(A^{l}A^{l+}+I)^{-1}\quad\mbox{і}\quad S=-(A^{l}A^{l+}+I)^{-1}.
$$
Використовуючи \eqref{6f41}, одержимо обмежені оператори
\begin{gather*}
R:Z_{0}=H^{-2ql}(\Omega)\to H^{2ql-2ql-2q(l-1)}(\Omega)=X_{0},\\
R:Z_{1}=H^{\lambda}(\Omega)\to H^{\lambda+4ql-2ql-2q(l-1)}(\Omega)=X_{1},\\
S:Z_{0}=H^{-2ql}(\Omega)\to H^{2ql}(\Omega)\hookrightarrow
H^{\lambda}(\Omega)=Y_{0},\\
S:Z_{1}=H^{\lambda}(\Omega)\to H^{\lambda+4ql}(\Omega)\hookrightarrow H^{\lambda}(\Omega)=Y_{1};
\end{gather*}
тут вкладення неперервні. Окрім того,
$$
AR=AA^{l-1}A^{l+}(A^{l}A^{l+}+I)^{-1}=
(A^{l}A^{l+}+I-I)(A^{l}A^{l+}+I)^{-1}=I+S
$$
на просторі $E=H^{-2ql}(\Omega)$. Таким чином, для введених операторів $R$ і $S$ виконуються умови (v)\,--\,(vii) твердження~2.

Згідно з цим твердженням і на підставі \eqref{6f39}, отримаємо рівності
\begin{equation*}
\bigl[H_{A,\lambda}^{-2q(l-1)}(\Omega),H^{\lambda+2q}(\Omega)\bigr]_{\psi}=
\bigl[(X_{0})_{T,Y_{0}},(X_{1})_{T,Y_{1}}\bigr]_{\psi}=
\bigl([X_{0},X_{1}]_{\psi}\bigr)_{T,[Y_{0},Y_{1}]_{\psi}}.
\end{equation*}
Тут на підставі теореми 3.2 (про інтерполяцію  з функціональним параметром соболєвських просторів на $\Omega$), наведеної у \cite{MikhailetsMurach14}, маємо рівність просторів
\begin{equation*}
[X_{0},X_{1}]_{\psi}=
\bigl[H^{-2q(l-1)}(\Omega),H^{\lambda+2q}(\Omega)\bigr]_{\psi}=
\bigl[H^{s-\varepsilon}(\Omega),H^{s+\delta}(\Omega)\bigr]_{\psi}=
H^{s,\varphi}(\Omega).
\end{equation*}
Ці рівності виконуються з точністю до еквівалентності норм. Окрім того,
$$
[Y_{0},Y_{1}]_{\psi}=
\bigl[H^{\lambda}(\Omega),H^{\lambda}(\Omega)\bigr]_{\psi}=
H^{\lambda}(\Omega).
$$
З останніх трьох виносних формул негайно випливає потрібна рівність \eqref{6f38}.

Теорему~1 доведено.

\textbf{\textit{Доведення теореми} 2.}
Виберемо довільно функцію $\chi\in C^{\infty}(\overline{\Omega})$ таку, що $\mathrm{supp}\,\chi\subset\Omega_0\cup\Gamma_0$. Виберемо також  функцію $\eta\in C^{\infty}(\overline{\Omega})$, яка задовольняє умови $\mathrm{supp}\,\eta\subset\Omega_0\cup\Gamma_0$ і $\eta=1$ у деякому околі $V$ множини $\mathrm{supp}\,\chi$ (звісно, цей окіл розглядається у топології на $\overline{\Omega}$). За умовою, $u\in H^{\sigma}_{A,\lambda}(\Omega)$ для деяких цілого числа $\sigma<\min\{s,2q\}$ і дійсного числа $\lambda>m+1/2-2q$ та, крім того,
\begin{equation}\label{proof2-a}
\eta(f,g)\in\left\{
\begin{array}{ll}
\mathcal{H}_{\lambda,s,\varphi}(\Omega,\Gamma),&
\hbox{якщо}\;\;s\leq m+1/2; \\
\mathcal{H}_{s-2q,s,\varphi}(\Omega,\Gamma), & \hbox{якщо}\;\;s>m+1/2.
\end{array}
\right.
\end{equation}
Тут, звісно, $g:=(g_{1},\ldots,g_{q})$ і $\eta(f,g):=(\eta f, (\eta\!\upharpoonright\!\Gamma)g_{1},\ldots,
(\eta\!\upharpoonright\!\Gamma)g_{q})$. У випадку, коли $s>m+1/2$, виберемо число $\lambda>m+1/2-2q$ так, щоб додатково виконувалась нерівність $\lambda<s-2q$.

Взагалі кажучи, $\chi u\notin H^{\sigma}_{A,\lambda}(\Omega)$; тому замість простору $H^{\sigma}_{A,\lambda}(\Omega)$ будемо використовувати більш широкій простір Соболєва\,--\,Ройтберга $H^{\sigma,(r)}(\Omega)$, замкнений відносно операції множення на довільну функцію класу $C^{\infty}(\overline{\Omega})$. Як було показано у доведенні теореми~1, норми у просторах $H^{\sigma}_{A,\lambda}(\Omega)$ і $H^{\sigma,(r)}_{A,\lambda}(\Omega)$ еквівалентні на щільному лінійному многовиді $C^{\infty}(\overline{\Omega})$. Отже, ці простори рівні з точністю до еквівалентності норм і тому виконується неперервне вкладення $H^{\sigma}_{A,\lambda}(\Omega)\hookrightarrow H^{\sigma,(r)}(\Omega)$.
Отож, нетерів оператор \eqref{6f26-m+1}, де беремо число $\sigma$ замість $s$, є розширенням нетерового оператора
\begin{equation}\label{proof2-c}
(A,B):H^{\sigma}_{A,\lambda}(\Omega)\to
\mathcal{H}_{\lambda,\sigma}(\Omega,\Gamma).
\end{equation}
Нагадаємо \cite{Roitberg96} (наслідок 2.3.1), що оператор множення на функцію класа $C^{\infty}(\overline{\Omega})$ є неперервним на кожному просторі Соболєва\,--\,Ройтберга.

Оскільки оператор \eqref{proof2-c} нетерів, а множина
$C^{\infty}(\overline{\Omega})\times(C^{\infty}(\Gamma))^{q}$ щільна в просторі $\mathcal{H}_{\lambda,\sigma}(\Omega,\Gamma)$, то за лемою Гохберга\,--\,Крейна \cite{HohbergKrein57} (лема 2.1) існує скінченновимірний простір $N_{0}\subset C^{\infty}(\overline{\Omega})\times(C^{\infty}(\Gamma))^{q}$ такий, що
\begin{equation*}
\mathcal{H}_{\lambda,\sigma}(\Omega,\Gamma)=
(A,B)\bigl(H^{\sigma}_{A,\lambda}(\Omega)\bigr)\dotplus N_{0}.
\end{equation*}
Позначимо через $P_{0}$ проектор простору $\mathcal{H}_{\lambda,\sigma}(\Omega,\Gamma)$ на перший доданок у цій прямій сумі паралельно другому доданку.

На підставі умови \eqref{proof2-a}, теореми~1 і твердження~1 робимо висновок, що у випадку $s\leq m+1/2$ виконується включення
$$
P_{0}(\eta(f,g))\in\mathcal{H}_{\lambda,s,\varphi}(\Omega,\Gamma)\cap
(A,B)\bigl(H^{\sigma}_{A,\lambda}(\Omega)\bigr)=
(A,B)\bigl(H^{s,\varphi}_{A,\lambda}(\Omega)\bigr),
$$
а у випадку $s>m+1/2$~--- включення
$$
P_{0}(\eta(f,g))\in\mathcal{H}_{s-2q,s,\varphi}(\Omega,\Gamma)\cap
(A,B)\bigl(H^{\sigma}_{A,\lambda}(\Omega)\bigr)=
(A,B)\bigl(H^{s,\varphi}(\Omega)\bigr).
$$
Тому $P_{0}(\eta(f,g))=(A,B)u_1$ для деякого $u_{1}\in H^{s,\varphi}(\Omega)\cap H^{\sigma}_{A,\lambda}(\Omega)$.
Тоді
$$
\eta(f,g)=P_{0}(\eta(f,g))+(I-P_{0})(\eta(f,g))=(A,B)u_1+u_{1}^{\circ},
$$
де $I$~--- тотожній оператор, а $u_{1}^{\circ}:=(I-P_{0})(\eta(f,g))\in N_{0}$. Крім того, за умовою, $(A,B)u=(f,g)$. Звідси отримуємо, що
$$
(A,B)(u-u_1)=(1-\eta)(f,g)+u_{1}^{\circ}.
$$
Оскільки $1-\eta=0$ на $V$ і $u_{1}^{\circ}\in N_{0}\subset C^{\infty}(\overline{\Omega})\times(C^{\infty}(\Gamma))^{q}$, а $\mathrm{supp}\,\chi\subset V$, то за теоремою про локальне підвищення регулярності розв'язків еліптичних крайових задач у просторах Соболєва\,--\,Ройтберга \cite{Roitberg96} (теорема 7.2.1) виконується включення
$$
w:=\chi(u-u_1)\in\bigcap_{\theta\geq r}H^{\theta,(r)}(\Omega)=
C^{\infty}(\overline{\Omega}).
$$
Таким чином, $\chi u=\chi u_1+w\in H^{s,\varphi}(\Omega)$, бо $u_1\in H^{s,\varphi}(\Omega)$. Отже, $u\in H^{s,\varphi}_{\mathrm{loc}}(\Omega_{0},\Gamma_{0})$ з огляду на довільність вибору функції $\chi$.

Теорему 2 доведено.


\begin{thebibliography}{99}

\bibitem{Venttsel59}
\emph{Вентцель A. Д.} О граничных условиях для многомерных диффузионных процессов~// Теория вероятн. и ее примен.~-- 1959.~-- \textbf{4}.~-- P. 172~-- 185. (Переклад англійською: \emph{Ventcel' A. D.} On boundary conditions for multi-dimensional diffusion processes~// Theory Probab. Appl.~-- 1959.~-- \textbf{4}.~-- P. 164~-- 177.)

\bibitem{Krasil'nikov61}
\emph{Красильников  В. Н.} О решении некоторых гранично-контактных задач линейной гидродинамики~// Прикл. мат. мех.~-- 1961.~-- \textbf{25}, № 4.~-- С. 764~-- 768. (Переклад англійською:
\emph{Krasil'nikov V. N.} On the solution of some boundary-contact problems of linear hydrodynamics~// J. Appl. Math. Mech.~-- 1961.~-- \textbf{25}, № 4.~-- P. 1134~-- 1141.)

\bibitem{VeshevKouzov77}
\emph{Вешев В.А., Коузов Д.П.} О влиянии среды на колебания пластин, сочлененных под прямым углом~// Акустический жуpнал.~-- 1977.~-- \textbf{23}, №~3.~-- С. 368~-- 377. (Переклад англійською:
\emph{Veshev V. A., Kouzov D. P.} Influence of the medium on the vibrations of plates joined at right angles~// Acoustical Physics.~-- 1977.~-- \textbf{23}, №~3.~-- P. 206~-- 211.)

\bibitem{Roitberg96}
\emph{Roitberg Ya.~A.}  Elliptic boundary value problems in the spaces of distributions.~-- Dordrecht: Kluwer Acad. Publisher, 1996.~-- xii+415~p.

\bibitem{KozlovMazyaRossmann97}
\emph{Kozlov V. A., Maz'ya V. G., Rossmann J.}  Elliptic boundary value problems in domains with point singularities.~-- Providence: Amer. Math. Soc., 1997.~-- 414~p.

\bibitem{Hermander63}
\emph{H\"ormander L.} Linear partial differential operators.~-- Berlin: Springer, 1963.~-- 285~p. (Переклад російською: \emph{Хермандер~Л.} Линейные дифференциальные операторы с частными производными.~-- Москва: Мир, 1965.~-- 380~с.)

\bibitem{Hermander83}
\emph{H\"ormander L.} The analysis of linear partial differential operators. II: Differential operators with constant coefficients.-- Berlin: Springer, 1983.~-- viii+391~p. (Переклад російською: \emph{Хермандер~Л.} Анализ линейных дифференциальных операторов с частными производными. Т.~2.~-- Москва: Мир, 1986.~-- 456~с.)

\bibitem{Jacob010205}
\emph{Jacob N.} Pseudodifferential operators and Markov processes: In 3
vol\-umes.~-- London: Imperial College Press, 2001, 2002, 2005.~-- xxii+493~p., xxii+453~p., xxviii+474~p.

\bibitem{MikhailetsMurach14}
\emph{Mikhailets V. A., Murach A. A.} H\"ormander spaces, interpolation, and elliptic problems.~-- Berlin, Boston: De Gruyter, 2014.~-- xii+297~p. (Видання російською доступне як arXiv:1106.3214.)

\bibitem{NicolaRodino10}
\emph{Nicola F., Rodino L.} Global Pseudodifferential Calculas on Euclidean spaces.~-- Basel: Birkh\"aser, 2010.~-- x+306~p.

\bibitem{Paneah00}
\emph{Paneah B.} The oblique derivative problem. The Poincar\'e problem.-- Berlin: Wiley--VCH, 2000.~-- 348~p.

\bibitem{Stepanets05}
\emph{Stepanets A. I.} Methods of approximation theory.~-- Utrecht: VSP, 2005.

\bibitem{Triebel01}
\emph{Triebel H.} The structure of functions.~-- Basel: Birkh\"aser, 2001.~-- xii+425~p.

\bibitem{MikhailetsMurach05UMJ5}
\emph{Mikhailets V. A., Murach A. A.} Elliptic operators in a refined scale of functional spaces~// Ukrainian Math. J.~-- 2005.~-- \textbf{57}, №~5.~-- P.~817~-- 825.

\bibitem{MikhailetsMurach06UMJ3}
\emph{Mikhailets V. A., Murach A. A.} Refined scales of spaces and elliptic boundary-value problems. II~// Ukrainian Math. J.~-- 2006.~-- \textbf{58}, №~3.~-- P. 398~-- 417.

\bibitem{MikhailetsMurach06UMJ11}
\emph{Mikhailets V. A., Murach A. A.} Regular elliptic boundary-value problem for homogeneous equation in two-sided refined scale of spaces~// Ukrainian Math.~J.~-- 2006.~-- \textbf{58}, №~11.~-- P. 1748~-- 1767.

\bibitem{MikhailetsMurach06UMB4}
\emph{Mikhailets V. A., Murach A. A.} Elliptic operator with homogeneous regular boundary conditions in two-sided refined scale of spaces~// Ukr. Math. Bull.~-- 2006.~-- \textbf{3}, №~4.~-- P. 529~-- 560.

\bibitem{MikhailetsMurach07UMJ5}
\emph{Mikhailets V. A., Murach A. A.} Refined scales of spaces and elliptic boundary-value problems. III~// Ukrainian Math. J.~-- 2007.~-- \textbf{59}, №~5.~-- P. 744~-- 765.

\bibitem{MikhailetsMurach08UMJ4}
\emph{Mikhailets V. A., Murach A. A.} An elliptic boundary-value problem in a two-sided refined scale of spaces.~-- Ukrainian Math. J.~-- 2008.~-- \textbf{60}, №~4.~-- P.~574~-- 597.

\bibitem{MikhailetsMurach12BJMA2}
\emph{Mikhailets V. A., Murach A. A.} The refined Sobolev scale,
inter\-po\-la\-tion, and elliptic problems~// Banach J. Math.
Anal.~-- 2012.~-- \textbf{6}, №~2.~-- P. 211~-- 281.

\bibitem{AnopMurach14MFAT2}
\emph{Anop A. V., Murach A. A.} Parameter-elliptic problems and interpolation with a function parameter~// Methods Funct. Anal. Topology.~-- 2014.~--- \textbf{20}, No~2.~-- P. 103--116.

\bibitem{AnopMurach14UMJ7}
\emph{Anop A. V., Murach A. A.} Regular elliptic boundary-value problems in the extended Sobolev scale~// Ukrainian Math. J.~-- 2014.~-- \textbf{66}, №~7.~-- P. 969~-- 985.

\bibitem{AnopKasirenko16MFAT4}
\emph{Anop A. V., Kasirenko T. M.} Elliptic boundary-value problems in H\"ormander spaces~// Methods Funct. Anal. Topology.~-- 2016.~-- \textbf{22}, №~4.~-- P. 295~-- 310.

\bibitem{MikhailetsMurach13UMJ3}
\emph{Mikhailets V. A., Murach A. A.} Extended Sobolev scale and elliptic operators~// Ukrainian Math.~J.~-- 2013.~-- \textbf{65}, №~3.~-- P. 435~-- 447.

\bibitem{MikhailetsMurach15ResMath1}
\emph{Mikhailets~V.~A., Murach~A.~A.} Interpolation Hilbert spaces between Sobolev spaces~// Results Math.~-- 2015.~-- \textbf{67}, №~1.~-- P. 135~-- 152.

\bibitem{DenkFaierman17arxiv}
\emph{Denk R., Faierman M.} An elliptic boundary problem acting on generalized Sobolev spaces~// arXiv:1710.01959v1.~-- 2017.~-- 21~p.

\bibitem{KasirenkoMurach17UMJ11}
\emph{Касiренко Т. М., Мурач О. О.} Еліптичні задачі з крайовими умовами високих порядків у просторах Хермандера~// Укр. мат. журн.~-- 2017.~-- \textbf{69}, №~11.~-- С. 1486~-- 1504.

\bibitem{LionsMagenes62}
\emph{Lions J.-L., Magenes E.} Probl\`emes aux limites  non homog\'enes, V~// Ann. Scuola Norm. Sup. Pisa (3).~-- 1962.~-- \textbf{16}.~-- P. 1~-- 44. (Italian)

\bibitem{LionsMagenes63}
\emph{Lions J.-L., Magenes E.} Probl\`emes aux limites non homog\'enes, VI~// J. Anal. Math.~-- 1963.~-- \textbf{11}.~-- P. 165~-- 188.

\bibitem{LionsMagenes71}
\emph{Лионс Ж.-Л., Мадженес Э.} Неоднородные граничные задачи и их приложения.~-- Москва: Мир, 1971.~-- 372~с.

\bibitem{Roitberg68}
\emph{Ройтберг Я. А.} Теоремы о гомеоморфизмах, осуществляемых эллиптическими операторами~// Доклады АН СССР.~-- 1968.~-- \textbf{180}, №~3.~-- С. 542--545. (Переклад англійською: \emph{Roitberg Ja. A.} Theorems on homeomorphisms which can be realized by elliptic operators~// Dokl. Math.~-- 1968.~-- \textbf{9}.~-- P. 656--660).

\bibitem{KostarchukRoitberg73}
\emph{Костарчук Ю.В., Ройтберг Я.А.} Теореми про ізоморфізми для еліптичних граничних задач з граничними умовами, які не є нормальними~// Укр. мат. журн.~-- 1973.~-- \textbf{25}, №~2.~-- С.~271~-- 277. (Переклад англійською: \emph{Kostarchuk Ju.(Yu.) V., Roitberg Ja.(Ya.) A.} Isomorphism theorems for elliptic boundary value problems with boundary conditions that are not normal~// Ukrainian Math. J.~-- 1973.~-- \textbf{25}, №~2.~-- P. 222~-- 226.)

\bibitem{Murach09MFAT2}
\emph{Murach A. A.} Extension of some Lions--Magenes theorems~// Methods Funct. Anal. Topology.~-- 2009.~-- \textbf{15}, №~2.~-- P. 152~-- 167.

\bibitem{MikhailetsMurach11Dop4}
\emph{Михайлец В. А., Мурач А. А.} Индивидуальные теоремы о разрешимости эллиптических за-дач и пространства Хермандера~// Доповіді НАН України.~-- 2011.~-- № 4.~-- С. 30~-- 36.

\bibitem{MurachChepurukhina15UMJ}
\emph{Murach A. A., Chepurukhina I. S.} Elliptic boundary-value problems in the sense of Lawruk on Sobolev and H\"ormander spaces~// Ukrainian Math.~J.~-- 2015.~-- \textbf{67}, №~5.~-- P. 764~-- 784.

\bibitem{Agranovich97}
\emph{Agranovich M. S.} Elliptic boundary problems~// Encycl. Math. Sci. Vol. 79. Partial differential equations, IX.~-- Berlin: Springer, 1997.~-- P. 1~-- 144.

\bibitem{FunctionalAnalysis72}
\emph{Функциональный анализ}~/ Под общ. ред. С.~Г.~Крейна.~-- Москва:
Наука, 1972.~-- 544~с.

\bibitem{Seneta76}
\emph{Seneta E.} Regularly varying functions.~-- Berlin: Springer, 1976.~-- 112~p. (Переклад російською: \emph{Сенета~Е.} Правильно меняющиеся функции.~-- М.: Наука, 1985.~-- 144~с.)

\bibitem{BinghamGoldieTeugels89}
\emph{Bingham N. H., Goldie C.~M., Teugels J.~L.} Regular variation.~-- Cambridge: Cambridge Univ. Press, 1989.~-- 512~p.

\bibitem{VolevichPaneah65}
\emph{Волевич Л. Р., Панеях Б. П.} Некоторые пространства обобщенных
функций и теоремы вложения~// Успехи мат. наук.~-- 1965.~--
\textbf{20}, №~1.~-- С. 3~-- 74. (Переклад англійською:
\emph{Volevich L. R., Paneah B. P.} Certain spaces of generalized functions and embedding theorems~// Russian Math. Surveys.~-- 1965.~-- \textbf{20}, №~1.~-- P. 1~-- 73.)

\bibitem{Roitberg64}
\emph{Ройтберг Я. А.} Эллиптические задачи с неоднородными граничными условиями и локальное повышение гладкости вплоть до границы обобщенных решений~// Доклады АН СССР.~-- 1964.~-- \textbf{157}, №~4.~-- С. 798~-- 801.

\bibitem{HohbergKrein57}
\emph{Гохберг И. Ц., Крейн М. Г.} Основные положения о дефектных числах, корневых векторах и индексах линейных операторов~// Успехи матем. наук.~-- 1957.~-- \textbf{12}, №~2.~-- С. 43~-- 118. (Переклад англійською: \emph{Gohberg~I.~C., Krein~M.~G.} The basic propositions on defect numbers, root numbers, and indices of linear operators~// Amer. Math. Soc. Transl., Ser.~2.~-- 1960.~-- \textbf{13}.~-- P. 185~-- 264.)

\bibitem{MikhailetsMurach08MFAT1}
\emph{Mikhailets V. A., Murach A. A.} Interpolation with a function parameter and refined scale of spaces~// Methods Funct. Anal. Topology.~-- 2008.~-- \textbf{14}, №~1.~-- P. 81~-- 100.

\end{thebibliography}
\end{document}